\documentclass{amsart}
\RequirePackage[OT1]{fontenc}
\RequirePackage{amsthm,amsmath,amssymb}
\RequirePackage{graphicx}
\RequirePackage[numbers]{natbib}
\RequirePackage[colorlinks,citecolor=blue,urlcolor=blue]{hyperref}

\newtheorem{defi}{Definition}
\newtheorem{lem}{Lemma}
\newtheorem{theo}{Theorem}
\newtheorem{rem}{Remark}
\newtheorem{prop}{Proposition}

\newcommand{\abs}[1]{{\left| #1\right|}}
\newcommand{\entier}[1]{{\lfloor #1\rfloor}}
\newcommand{\norme}[1]{ {\left\lVert  #1\right\rVert}}

\newcommand{\dd}{\text{d}}
\newcommand{\Var}{\operatorname{Var}}
\newcommand{\Cov}{\operatorname{Cov}}

\newcommand{\HReg}{{\bf HReg}}
\newcommand{\HErg}{{\bf HErg}}
\newenvironment{pf}{\ \\ {\bf Proof }}{\hfill\mbox{$\square$}\medskip}

\newcommand{\modar}{}

\numberwithin{equation}{section}
\theoremstyle{plain}

\begin{document}

\title[Rate of estimation for stochastic damping Hamiltonian systems]{Rate of estimation for the stationary distribution of stochastic damping Hamiltonian systems with continuous observations
	}

\author{Sylvain Delattre}
\address{Laboratoire de Probabilit\'es Statistique et Mod\'elisation, UMR 7599, Universit\'e Denis Diderot, 77251 Paris Cedex
		05, France}
\email{}

\author{Arnaud Gloter}
\address{Laboratoire de Math\'ematiques et Mod\'elisation d'Evry, CNRS, Univ Evry, Universit\'e Paris-Saclay, 91037, Evry, France}
\curraddr{}
\email{}
\thanks{}

\author{Nakahiro Yoshida}
\address{Graduate School of Mathematical Sciences, University of Tokyo: 3-8-1 Komaba, Meguro-ku, Tokyo 153-8914, Japan.}
\curraddr{}
\email{}
\thanks{}



\date{26, january, 2020}

\maketitle

\begin{abstract}
		We study the problem of the non-parametric estimation for 
		the density $\pi$ of the stationary distribution of a stochastic two-dimensional damping Hamiltonian system $(Z_t)_{t\in[0,T]}=(X_t,Y_t)_{t \in [0,T]}$. From the continuous observation of the sampling path on $[0,T]$, we study the rate of estimation for $\pi(x_0,y_0)$ as $T \to \infty$. We show that kernel based estimators can achieve the rate $T^{-v}$ for some explicit exponent $v \in (0,1/2)$. One finding is that the rate of estimation depends on the smoothness of $\pi$ and is completely different with the rate appearing in the standard i.i.d.\ setting or in the case of two-dimensional non degenerate diffusion processes. Especially, this rate depends also on $y_0$. Moreover, we obtain a minimax lower bound on the $L^2$-risk for pointwise estimation, with the same rate $T^{-v}$, up to $\log(T)$ terms.
		
\end{abstract}

\section{Introduction}
The class of 
hypo-elliptic diffusion processes, for which the diffusion coefficient is degenerate,
has been the subject of many recent works and 
is used for modeling
in many fields, such as mathematical finance, biology, neuro-science,
mechanics,  ecology,...
 (see e.g. \cite{Ditlevsen_Samson_19}, \cite{Leon_Samson18}, \cite{Ditlevsen_Sorensen04} and references therein).
  In this paper, we focus on the situation of a bi-dimensional hypo-elliptic process, describing the evolution in time of the couple position/velocity of some quantity.
The velocity $Y$ is modeled by a non degenerate one-dimensional diffusion process, while the position $X$ is its integral, and the resulting bi-dimensional process $Z=(X,Y)$ is hypo-elliptic.  
Our aim is to estimate the density $\pi$ of the stationary measure of this diffusion process, under the assumption of an ergodic setting.

The problem of non-parametric estimation of the stationary measure of a continuous mixing
process is a long-standing problem (see for instance N'Guyen \cite{NGuyen79}, or Comte and Merlevede \cite{ComMer02} and references therein). Based on sample of length $T$ of the data
Comte and Merlevede \cite{ComMer02} find estimators of the stationary measure, 
 converging at rate depending on the smoothness of the stationary measure and slower than $\sqrt{T}$, as it is usual in non-parametric problems.

In the specific context where the continuous time process is a one-dimensional diffusion process, observed continuously on some interval $[0,T]$, the finding is different.  It is shown that the rate of estimation of the stationary measure is $\sqrt{T}$ 
(see  Kutoyants \cite{Kutoyants04}). 
The rate of estimation is thus independent of the smoothness of the object that one estimates, in contrast to the typical non-parametric situation. The optimal estimator is very specific to the diffusive nature of the process as it relies on the local time of the process. Remark that if the process is a diffusion observed discretely on $[0,T]$ with a sufficiently high frequency it is possible to estimate with rate $\sqrt{T}$ also (see \cite{Nishiyama11}, \cite{ComMer05}).

The case of multi-dimensional non-degenerate diffusions is treated in Dalalyan and Rei\ss{} \cite{DalRei07} and Srauss \cite{Strauch18}. In that case, the local time process is not available, but it is shown that for non degenerate diffusion process of dimension $d=2$, there exists an estimator of the pointwise values of the stationary measure with rate $\sqrt{T}/\log(T)^2$. This rate of estimation does not depend on the smoothness of the stationary measure for $d=2$.
   In the situation of a non degenerate diffusion with dimension $d \ge 3$,  
they find estimators whose rate is polynomial in $T$, depends on both the smoothness of the stationary measure and the dimension $d \ge 3$. For $d \ge 3$, the rate is strictly slower than $\sqrt{T}$, but faster than the rate appearing in standard multivariate density estimation from $T$ i.i.d. observations. Hence, for $d \ge 3$, the diffusive structure of the process, enables to get a faster estimation rate than for the i.i.d. case, as well.

In the case of hypo-elliptic processes, fewer results exist for the estimation of the stationary distribution.  In \cite{Comte_et_al_17}, the authors consider the case of two-dimensional process $Z=(X,Y)$, where $Y$ is a velocity and $X$ a position. Based on a discrete sampling of the path of size $n$, they propose an estimator of the stationary measure which converges at a non parametric rate depending on the smoothness of the stationary measure. Assuming that the stationary measure has anisotropic regularity $(k_1,k_2)$, the proposed estimator has a rate depending on the harmonic mean $2/(1/k_1+1/k_2)$ as it is for the optimal rate of estimation of distribution on in the case of i.i.d. sequence.  

In this paper, we focus on the situation where the process is observed continuously on $[0,T]$ and our goal is to determine what is the optimal rate of estimation of $\pi$ in this context.
Assuming that the stationary density $(x,y) \mapsto \pi(x,y)$ has 
an anisotropic 
H\"older regularity with index $k_1$ with respect to the variable $x$, and $k_2$ 
with respect to $y$,
we construct an estimator of $\pi(x_0,y_0)$ based on $(X_t,Y_t)_{t \in [0,T]}$.
This estimator achieves some rate $T^{-v(k_1,k_2)}$ and this rate of convergence 
depends on the smoothness of $\pi$ in a very specific way. 
Indeed, the expression of the rate of estimation involves only $k_1$ or $k_2$, depending on the relative positions of these two smoothness indexes. 
This shows the specificity of the estimation problem for continuous observation of process $Z=(X,Y)$ with a non degenerate diffusive velocity $Y$, and a degenerate component $X$.
It is noteworthy that we find a rate of estimation 
slower than in the case of continuous observation of a non-degenerate diffusion in dimension $2$.   Another interesting finding is that the rate of estimation of $\pi(x_0,y_0)$ depends on the point where one estimates the stationary measure, and is slower for points corresponding to null velocity $y_0=0$.
A crucial ingredient in the study of the rate of estimation is to derive the variance of a kernel estimator 
with a choice of bandwidth $(h_1,h_2)$. We show that the variance of the kernel estimator depends in a completely unsymmetrical way on $h_1$ and $h_2$. As a consequence, we get that the optimal bandwidth choice for the estimator is such that one of the two bandwidth can go almost arbitrarily fast to $0$, while the optimal choice for the other bandwidth depends sharply on the smoothness of the stationary measure. 

Also, we show a lower bound for the minimax risk of estimation of $\pi(x_0,y_0)$ on a class of  hypo-elliptic diffusion models with stationary measure of H\"older regularity  $(k_1,k_2)$. This proves that it is impossible to estimate uniformly on this  class with a rate faster than $T^{-v(k_1,k_2)}$ (up to $\log(T)$ term).

The outline of the paper is the following. In Section \ref{S:Assumptions}, we present the model and give assumptions that are sufficient to get an ergodic system with stationary measure admitting a  density $\pi$. In Section \ref{S:Estimator}, we present the construction of the estimator and states the results on their rate of convergence (in Theorem \ref{T:upper_bound} for $y_0 \neq 0$, and
Theorem \ref{T:upper_bound_y0_0} for $y_0 =0$). In Section \ref{S:Variance}, we prove the upper bound on the variance of the kernel estimator. We also illustrate the very specific behaviour of the variance of the kernel estimator for hypo-elliptic diffusion  by numerical simulations. 
In Section \ref{S:minimax}, we state and prove minimax lower bounds for the risk of estimation.
In the Appendix, we prove some technical results 
used  in the proofs of Section \ref{S:Variance}.

%
%
%
%
%
%

\section{{\modar Hamiltonian system and mixing property}}

\label{S:Assumptions}
Let us consider $(\Omega,\mathcal{A},\mathbb{P})$, some probability space on which a standard one dimensional Brownian motion $(B_t)_{t\ge0}$ is defined. 
We assume that the process $(Z_t)_{t \ge 0}=(X_t,Y_t)_{t \ge 0}$ is solution of the stochastic differential equation
\begin{align} \label{E:model_X}
&\dd X_t = Y_t \dd t \\
\label{E:model_Y}
&\dd Y_t = a(X_t,Y_t) \dd B_t - [\beta(X_t,Y_t) Y_t + V'(X_t)] \dd t,
\end{align}
where $(X_0,Y_0)$ is a random variable independent of $(B_t)_t$.

We introduce the following regularity assumption on the coefficients.

{\modar
{Assumption \HReg :} 
\begin{itemize}
	\item
The functions $a : \mathbb{R}^2 \to \mathbb{R}$ is a $\mathcal{C}^\infty$ function and $\overline{a}>a(x,y)\ge \underline{a} >0$, 
 and for some constants $\overline{a}>\underline{a}>0$.

\item The function 
$\beta :  \mathbb{R}^2 \to \mathbb{R}$ is continuously differentiable, and such that
$\abs{\beta(x,y)} \le \overline{\beta}$ for all $(x,y) \in \mathbb{R}^2$, 
and for some constant
$\overline{\beta} >0$.
Moreover, we have
$\beta(x,y) > \underline{\beta}$, $\forall x \in [l,\infty),y\in \mathbb{R}$, where $l\ge 0$, $\underline{\beta}>0$ are two constants.

\item The function $V: \mathbb{R}\to \mathbb{R}$ is lower bounded, and with $\mathcal{C}^2$ regularity.
\end{itemize}
}

It is shown in \cite{Wu01}, that under the assumptions \HReg{} the S.D.E. \eqref{E:model_X}--\eqref{E:model_Y} admits a weak solution, which satisfies the Markov property, {\modar and the associated semi group $(P_t)_{t\ge 0}$ is strongly Feller \cite{Wu01}, which in turn implies that the process is strongly Markovian. Let us stress that the sign condition on $\beta$ for large $x$ together with the existence of lower bound on $V$ are crucial to insure that the solution of \eqref{E:model_X}--\eqref{E:model_Y} does not explode in finite time.} 
Of course, if we know that $(x,y)\mapsto a(x,y)$, $(x,y)\mapsto y\beta(x,y)$ and $x \mapsto V'(x)$ are globally Lipschitz, the solutions of the S.D.E. exists in the strong sense.

We now introduce an assumption on the potential $V$ of the system that ensure that the process tends to some equilibrium.


{Assumption \HErg :} one has $ \displaystyle \lim_{\abs{x}\to \infty} V'(x) \text{sign}(x) =+ \infty$.

Is is shown in \cite{Wu01} that under \HReg{} and \HErg, one can construct a Lyapounov function $\Psi \ge 1$, and that a stationary probability $\pi$ exists and is unique for the process $Z=(X,Y)$, and satisfies $\pi(\Psi)<\infty$. 
It is shown 
in \cite{Wu01} (see Theorem 2.4) that for some $D>0$ and $\rho >0$,
\begin{equation}\label{E:mixing_Wu}
\forall t \ge 0, \quad \forall z \in \mathbb{R}^2,  \quad \sup_{\abs{f} \le \Psi}  \abs{P_t(f)(z) - \int_{\mathbb{R}^2} f(z') \pi(\dd z')}
\le D \Psi(z) e^{-\rho t}
\end{equation} 
for any function measurable function $f$ such as $f/\Psi$ is bounded on $\mathbb{R}$ and where $(P_t)_{t \ge 0}$ is the semi group of $Z$,
$$
P_t(f)(z)=E[f(X_t,Y_t) \mid (X_0,Y_0)=z].
$$
Remark that under \HReg{} and \HErg, it is possible to construct a Lyapounov function such that $\Psi(x,y) \ge \frac{1}{C}
\exp \left(\frac{1}{C} [\abs{y}^2 + V(x) ] \right) $, for some constant $C>0$ (see (3.10) in \cite{Wu01}). As a consequence, \eqref{E:mixing_Wu} applies to functions $f$ with exponential growth, and the invariant distribution admits finite exponential moments.

{\modar Hence, we can state the following proposition}
\begin{prop} \label{P:HErg_HMixv2}
Assume that the coefficients of the equation \eqref{E:model_Y} satisfy \HReg{} and \HErg, then there exists a stationary solution $Z=(X,Y)$ to the S.D.E.  \eqref{E:model_X}--\eqref{E:model_Y}, and the stationary distribution is unique and admits some density $\pi$.
 Moreover, there exist constants $D_{\text{erg}}>0$ and $\rho>0$ such that for 
any  bounded measurable functions $f$, $g$, we have
\begin{equation}\label{E:Mixing}
\forall t \ge 0,\quad
\abs{\text{cov}(f(Z_0),g(Z_t))} \le {\modar D_{\text{erg}}} \norme{f}_\infty \norme{g}_\infty e^{-\rho t}.
\end{equation}
\end{prop}	
The Proposition \eqref{P:HErg_HMixv2} is a consequence of the results shown in \cite{Wu01}. Indeed, from the fact that the Lyapounov function $\Psi$ is greater than $1$ and integrable with respect to the stationary measure, one can check that \eqref{E:Mixing} is a consequence of \eqref{E:mixing_Wu}.

\section{Estimator and upper bounds}\label{S:Estimator}
In this section we introduce the expression for our estimator of the stationary measure $\pi$ of the S.D.E. \eqref{E:model_X}--\eqref{E:model_Y} and prove that the estimator achieves some rate of convergence, depending on the smoothness of $\pi$. 

Let $\varphi : \mathbb{R} \to \mathbb{R}$ a bounded, compactly supported function. {\modar For convenience, we suppose that the support of $\varphi$ is $[-1,1]$.}  We assume that 
\begin{equation}\label{E:oscillation_kernel}
\int_{\mathbb{R}} \varphi(u) \dd u=1, \quad \int_{\mathbb{R}} \varphi (u) u^l \dd u=0,
\text{ for $l \in \{1,\dots, L\}$ where $L\ge 1$.}
\end{equation} 

We let $h_1(T)>0$, $h_2(T)>0$ be two bandwidths which converge to zero as $T\to \infty$,
and we consider a kernel estimator of $\pi$ at the point $(x_0,y_0) \in \mathbb{R}^2$ as
\begin{equation} \label{E:Def_pi_hat}
\hat{\pi}_T(x_0,y_0) = \frac{1}{T} \int_0^T \varphi_{h_1(T),h_2(T)}(X_s-x_0,Y_s-y_0)\dd s,
\end{equation}
where 
\begin{equation}
\label{E:def_var_phi_double}
\varphi_{h_1,h_2}(x-x_0,y-y_0)=\frac{1}{h_1h_2} \varphi(\frac{x-x_0}{h_1})\varphi(\frac{y-y_0}{h_2}).
\end{equation}
We assume that the two bandwidths satisfy,
\begin{align} \label{E:h_sub_poly_1}
\exists K >0, \quad & h_1(T)^{-1}+h_2(T)^{-1} {\modar \le K(1+T^K),~ \forall T>0,}
\\  \label{E:h_sub_poly_2}
& \sqrt{h_1(T)}+h_2(T) {\modar \le K (\log(T))^{-3/2} \wedge 1,~ \forall T>1.}
\end{align}
The two previous conditions {\modar insure} that the bandwidths go faster to zero than {\modar the logarithmic rate} by \eqref{E:h_sub_poly_2}, but not faster than any polynomial rates by \eqref{E:h_sub_poly_1}. Actually these two bandwidths will be specified later (see equations \eqref{E:def_bandwidth_case_1_1}, \eqref{E:def_bandwidth_case_1_2}, \eqref{E:def_bandwidth_case_2_1}, \eqref{E:def_bandwidth_case_2_2} in the proofs of Theorems \ref{T:upper_bound} and \ref{T:upper_bound_y0_0}). 

We introduce the class of H{\"o}lder functions.
\begin{defi}
 For $(k_1,k_2) \in (0,\infty)^2$, and $R>0$, we denote $\mathcal{H}^{k_1,k_2}(R)$ the set of functions $f : \mathbb{R}^2 \to \mathbb{R}$ such that $x \mapsto f(x,y)$ and $y \mapsto f(x,y)$ are respectively of class $\mathcal{C}^{\entier{k_1}}$ and $\mathcal{C}^{\entier{k_2}}$, and satisfy the control, $\forall x,y$ in $\mathbb{R}$ and $h \in [-1,1]$,
\begin{align*} 
&\abs{f(x,y)} \le R,\\
&\abs{\frac{\partial^{\entier{k_1}} f}{\partial x^{\entier{k_1}}}(x+h,y)-\frac{\partial^{\entier{k_1}} f}{\partial x^{\entier{k_1}}}(x,y)} \le R \abs{h}^{k_1-\entier{k_1}}, \\
&\abs{\frac{\partial^{\entier{k_2}} f}{\partial x^{\entier{k_2}}}(x,y+h)-\frac{\partial^{\entier{k_2}} f}{\partial x^{\entier{k_2}}}(x,y)} \le R \abs{h}^{k_2-\entier{k_2}}.
\end{align*}
\end{defi}

We can state the main results on the asymptotic behaviour of the estimator. This behaviour 
is different according to the fact that we estimate the value of the stationary measure on a point $(x_0,y_0)$ corresponding to a null velocity or not.

\begin{theo} \label{T:upper_bound}
Assume that $Z=(X,Y)$ is a stationary solution to \eqref{E:model_X}--\eqref{E:model_Y} and that Assumptions \HReg, \HErg{} hold true. We assume that the stationary distribution $\pi$ belongs to the set $\mathcal{H}^{k_1,k_2}(R)$
for $k_1>0,k_2>0,R>0$, with $\max(k_1,k_2) \le L$ (recall \eqref{E:oscillation_kernel}).

Assume that $y_0 \neq 0$.
{\modar Then, there exist bandwidths $(h_1(T))_{T}$, $(h_2(T))_{T}$, depending only on $k_1$ and $k_2$, such that the estimator satisfies :
	\begin{align}
	\label{E:upper_bound_case_1_1}
	&\text{if $k_1 < k_2/2$,\quad}
	E\left[ (\hat{\pi}_T(x_0,y_0)-\pi(x_0,y_0))^2 \right] \le C 
	T^{-\frac{2 k_2}{2k_2+1}},
	\\
	\label{E:upper_bound_case_1_2}
	&\text{if $k_1 \geq k_2/2$,\quad}
	E\left[ (\hat{\pi}_T(x_0,y_0)-\pi(x_0,y_0))^2 \right] \le C 
	T^{-\frac{2 k_1}{2k_1+1/2}},
	\end{align}
for some constant $C$ independent of $T$.

}
\end{theo}

\begin{theo} \label{T:upper_bound_y0_0}
	Assume that $Z=(X,Y)$ is a stationary solution to \eqref{E:model_X}--\eqref{E:model_Y} and that Assumptions \HReg, \HErg{} hold true. We assume that the stationary distribution $\pi$ belongs to the set $\mathcal{H}^{k_1,k_2}(R)$
	for $k_1>0,k_2>0,R>0$, with $\max(k_1,k_2) \le L$ (recall \eqref{E:oscillation_kernel}).
%
	
	Assume that $y_0 = 0$.
	{\modar Then, there exist bandwidths $(h_1(T))_{T}$, $(h_2(T))_{T}$, depending only on $k_1$ and $k_2$, such that the estimator satisfies :
		\begin{align}
		 \label{E:upper_bound_case_2_1}
		 &\text{if $k_1 < k_2/3$,\quad}
		E\left[ (\hat{\pi}_T(x_0,y_0)-\pi(x_0,y_0))^2 \right] \le C 
		( \frac{T}{\log T})^{-\frac{2 k_2}{2k_2+2}},
		\\
		&\text{if $k_1 \geq k_2/3$,\quad} 
		\label{E:upper_bound_case_2_2}
		E\left[ (\hat{\pi}_T(x_0,y_0)-\pi(x_0,y_0))^2 \right] \le C 
		T^{-\frac{2 k_1}{2k_1+2/3}},
		\end{align}
		for some constant $C$ independent of $T$.
}	
\end{theo}
{\modar 
\begin{rem}
	The rates of estimation obtained in Theorems \ref{T:upper_bound}--\ref{T:upper_bound_y0_0} are completely different with the usual one in several ways. First, they do not depend on the harmonic mean of the smoothness index $k_1$, $k_2$, as it is usual in non-parametric setting. Second, the rate depends on the point $(x_0,y_0)$ where the density is estimated. We state in Section \ref{S:minimax} a minimax lower bound for the $L^2$ risk of estimation of $\pi(x_0,y_0)$ with the same rates (up to $\log$ terms).
\end{rem}
}
{\modar The asymptotic behaviour of the estimator relies on the standard bias variance decomposition. Hence, we need sharp evaluations for the variance of the estimator, that are stated below, and will be proved in Section \ref{S:Variance}.}

\begin{prop} \label{P:variance_y_0_non_zero}
	Assume that $Z=(X,Y)$ is solution to \eqref{E:model_X}--\eqref{E:model_Y}, that Assumptions \HReg, \HErg,
	hold true and {\modar $\norme{\pi}_\infty \le R$ for some $R>0$.}
	
	Assume that $y_0 \neq 0$.
	Then, there exists some constant $C$, such that for {\modar all $T>0$,}
	\begin{equation*}
	\Var( \hat{\pi}_T(x_0,y_0) )
	\le
C \left( \frac{1}{T h_2(T)} \wedge  \frac{1}{T \sqrt{h_1(T)}} \right) + C \varepsilon(T,h_1(T),h_2(T)),		
	\end{equation*}
where 
\begin{equation} \label{E:def_varepsilon_y0_non_zero}
\varepsilon(T,h_1(T),h_2(T)) \leq 
\frac{\abs{\log (h_1(T) h_2(T))}^{C}}{T}.
\end{equation}
\end{prop}
\begin{prop} \label{P:variance_y_0_zero}
	Assume that $Z=(X,Y)$ is a solution to \eqref{E:model_X}--\eqref{E:model_Y}, that Assumptions \HReg, \HErg, 
	hold true and {\modar that $\norme{\pi}_\infty \le R$ for some $R>0$.}
	
	Assume that $y_0 = 0$.
	Then, there exists some constant $C$, such that for all $T>0$,
	\begin{equation*}
	\Var( \hat{\pi}_T(x_0,y_0) )
	\le
	C \left(  \frac{\ln (T)}{T h_2(T)^2}\wedge \frac{1}{T h_1(T)^{2/3}} \right) + C \varepsilon(T,h_1(T),h_2(T)),		
	\end{equation*}
	where 
	\begin{equation*}
	\varepsilon(T,h_1(T),h_2(T)) \leq 
	\frac{\abs{\log (h_1(T) h_2(T) )}^{C}}{T}.
	\end{equation*}
\end{prop}	

We can now prove the that our estimator achieves the rates given in Theorems \ref{T:upper_bound}--\ref{T:upper_bound_y0_0}.
\begin{pf}[Proof of Theorem \ref{T:upper_bound}.]
	We write the usual bias-variance decomposition,
	\begin{multline} \label{E:bias_variance_decomp}
	E[(\hat{\pi}_T(x_0,y_0)- \pi(x_0,y_0))^2] \le \\
	\abs{E(\hat{\pi}_T(x_0,y_0))- \pi(x_0,y_0)}^2 + 
	\Var(\hat{\pi}_T(x_0,y_0)).
	\end{multline}	
	Using the stationarity of the process, 
	we can upper bound the bias term as 
	\begin{multline*}
	\abs{E(\hat{\pi}_T(x_0,y_0))- \pi(x_0,y_0)}^2   
	\\
	= \left( \int_{\mathbb{R}^2} \varphi(u)\varphi(v) 
	[
	\pi(x_0+u h_1(T),y_0+v h_2(T)) - \pi(x_0,y_0)] 
	\dd u \dd v
	\right)
	\\ 
	\le  C (h_1(T)^{2k_1}+h_2(T)^{2k_2}),
	\end{multline*}
	where in the last line we used $\pi \in \mathcal{H}^{k_1,k_2}(R)$ with  \eqref{E:oscillation_kernel} (see e.g. \cite{Tsybakov_book}, or Proposition 1 in \cite{Comte_Lacour13} for details). 
	{\modar 
		We now use 
		the results of Proposition \ref{P:variance_y_0_non_zero} on the variance of the estimator and
		choose the optimal bandwidths $h_1(T)$, $h_2(T)$.}

	$\bullet$ Case 1: $k_1<k_2/2$.
	
	Using Proposition \ref{P:variance_y_0_non_zero}
	\begin{multline*}
	E[(\hat{\pi}_T(x_0,y_0)- \pi(x_0,y_0))^2]  \le  C (h_1(T)^{2k_1}+h_2(T)^{2k_2}) \\+ \frac{C}{T h_2(T)} 
	+ C \varepsilon(T,h_1(T),h_2(T)). 
	\end{multline*}
	We now choose to balance $h_2(T)^{2k_2}$ with the main contribution of the variance term and let the contribution of $h_1(T)$ on the bias  be smaller. It yields us to set
	\begin{equation} \label{E:def_bandwidth_case_1_1}
	h_2(T)=T^{-1/(2k_2+1)}, \quad h_1(T)=T^{-C_1}, \text{ where $C_1 \ge {\frac{k_2}{k_1(2k_2+1)}}$}.
	\end{equation}
	With these choices and recalling \eqref{E:def_varepsilon_y0_non_zero}, we get \eqref{E:upper_bound_case_1_1}.
	
	$\bullet$ Case 2: $k_1\ge k_2/2$.
	
	We use again Proposition \ref{P:variance_y_0_non_zero}:
	\begin{multline*}
	E[(\hat{\pi}_T(x_0,y_0)- \pi(x_0,y_0))^2]  \le  C (h_1(T)^{2k_1}+h_2(T)^{2k_2})
	\\
	+ \frac{C}{T \sqrt{h_1(T)}} + C \varepsilon(T,h_1(T),h_2(T)). 
	\end{multline*}
	Balancing  the variance and bias terms yields to
	\begin{equation} \label{E:def_bandwidth_case_1_2}
	h_1(T)=T^{-1/(2k_1+1/2)}, \quad h_2(T)=T^{-C_2}, \text{ where $C_2 \ge {\frac{k_1}{k_2(2k_1+1/2)}}$},
	\end{equation}
	and \eqref{E:upper_bound_case_1_2} follows.
\end{pf} 

\begin{pf}[Proof of Theorem \ref{T:upper_bound_y0_0}]
	We use again the bias/variance decomposition \eqref{E:bias_variance_decomp} and exploit now the results of Proposition \ref{P:variance_y_0_zero}. 
	
	$\bullet$ Case 1: $k_1<k_2/3$.
	
	We have that,
	\begin{multline*}
	E[(\hat{\pi}_T(x_0,y_0)- \pi(x_0,y_0))^2]  \le  C (h_1(T)^{2k_1}+h_2(T)^{2k_2}) 
	\\
	+ \frac{C\ln(T)}{T h_2(T)^2} + C \varepsilon(T,h_1(T),h_2(T)). 
	\end{multline*}
	We set
	\begin{equation} \label{E:def_bandwidth_case_2_1}
	h_2(T)=(\frac{T}{\ln(T)})^{-1/(2k_2+2)}, \quad h_1(T)=T^{-C_1}, \text{ where $C_1 \ge  {\frac{k_2}{k_1(2k_2+2)}}$},
	\end{equation}
	and \eqref{E:upper_bound_case_2_1} follows.
	
	$\bullet$ Case 2: $k_1\ge k_2/3$.
	
	We have that,
	\begin{multline*}
	E[(\hat{\pi}_T(x_0,y_0)- \pi(x_0,y_0))^2]  \le  C (h_1(T)^{2k_1}+h_2(T)^{2k_2}) 
	\\ + \frac{C}{T h_1(T)^{2/3}} + C \varepsilon(T,h_1(T),h_2(T)). 
	\end{multline*}
	The choice
	\begin{equation} \label{E:def_bandwidth_case_2_2}
	h_1(T)=T^{-1/(2k_1+2/3)}, \quad h_2(T)=T^{-C_2}, \text{ where $C_2 \ge {\frac{k_1}{k_2(2k_1+2/3)}}$}
	\end{equation}
	gives \eqref{E:upper_bound_case_2_2}.
\end{pf} 
\begin{rem} 
	\begin{itemize}
		\item The optimal choices for the bandwidths are given in \eqref{E:def_bandwidth_case_1_1}, \eqref{E:def_bandwidth_case_1_2}, \eqref{E:def_bandwidth_case_2_1} and \eqref{E:def_bandwidth_case_2_2}. 
		In all the situations, 
		we see that one of the two bandwidths $h_1(T)$ or $h_2(T)$ can be chosen ``arbitrarily small'' (as the constants $C_1$ and $C_2$ in \eqref{E:def_bandwidth_case_1_1}, \eqref{E:def_bandwidth_case_1_2}, \eqref{E:def_bandwidth_case_2_1}, \eqref{E:def_bandwidth_case_2_2} can be arbitrarily large). In means that the bias induced by the variation of $\pi$ along one of the two variables $x$ or $y$ can be arbitrarily reduced by the choice of a very thin bandwidth. It explains why the expression of the rate of estimation depends only on one index of smoothness $k_1$ or $k_2$.
		\item The fact that one of the two bandwidths can be chosen arbitrarily small is reminiscent to the situation of the estimation of the stationary measure $\pi(z)$ for a one dimensional diffusion process $(Z_t)_{t \in [0,T]}$ observed continuously. In that case, the efficient estimator is based on the local time of the process (see \cite{Kutoyants04}) and the rate is $\sqrt{T}$ independently of the smoothness of $\pi$.
		The use of local time is a way 
		to give a rigorous 
		analysis of the quantity
		$\tilde{\pi}_T(z)=\frac{1}{T}\int_0^T \delta_{\{z\}}(Z_s) \dd s$, where $\delta_{\{z\}}$ is the Dirac mass located at $z$. We see that $\tilde{\pi}_T$ is essentially a kernel estimator with bandwidth $h=0$, for which the bias is reduced to $0$.		
	\end{itemize}
\end{rem}	

\section{Variance of the kernel estimator}\label{S:Variance}
In this section we prove the crucial upper bounds given in Propositions \ref{P:variance_y_0_non_zero}--\ref{P:variance_y_0_zero}.  We first need to state two lemmas related to the behaviour of density and semi group of the process, and whose proofs
	are postponed to Section \ref{Appendix}.
	\begin{lem}[Corollary 2.12 in \cite{Cattiaux_et_al_14}]
		\label{L:HShort}
		Assume \HReg{}. Then, the process admits a transition density $(p_t)_{t >0}$ which satisfies the following upper bound. For all $K$ compact subset of
		$\mathbb{R}^2$, $\forall (x,y,x',y') \in K\times K$,
		$\forall t \in (0,1)$
		\begin{equation} \label{E:Hshort}
		p_t((x,y);(x',y')) \le p_t^{G}((x,y);(x',y'))  +
		p_t^U((x,y);  (x' ,y'))
		\end{equation}
		where
		\begin{equation} \label{E:Def_p_G}
		p_t^{G}((x,y);(x',y')) =
		\frac{{\modar C_G}}{t^2} \exp \left( 
		- \frac{1}{{\modar C_G}} \left[ \frac{(y-y')^2}{t} + \frac{(x'-x-\frac{y+y'}{2}t)^2}{t^3} \right]
		\right),
		\end{equation}
		for some ${\modar C_G}>0$ and
		$p^U$ is a measurable non negative function such that for any compact $K \subset \mathbb{R}^2$ and $(x,y) \in  K$ we have for all $t \in (0,1)$,
		\begin{equation} \label{E:pte_p_U}
		\int_{\mathbb{R}^2} p_t^{U}((x,y);( x', y')) \dd x' \dd y' \leq C_U \exp( -t^{-1} C_U^{-1}),
		\end{equation}
		for $C_U>0$. The two constants $C_G$ and $C_U$ are independent of $t \in [0,1]$, but depend on the compact set $K$.
	\end{lem}
{\modar The Lemma \ref{L:HShort} gives us a control on the short time behaviour of the transition density.
	For the sequel, we need a control
	valid for any time. This is the purpose of the following lemma about the semi group of the process.
	\begin{lem} \label{L:HGlobal}
		Assume \HReg{}, and let $K$ be a compact subset of $\mathbb{R}^2$. Then, there exists a constant $\widetilde{C}_K$ such that for all $0<D<1$, $t\ge D$, $\forall z \in \mathbb{R}^2$, and any $f$ measurable bounded function with support on $K$,
		\begin{equation} \label{E:equa_type_HGlob}
		\abs{P_t(f)(z)} \le 
		\widetilde{C}_K \left[ \frac{\norme{f}_{{L}^1(\mathbb{R}^2)}}{D^2} + \norme{f}_{\infty} e^{-1/(\widetilde{C}_K D)} \right].
		\end{equation}
	\end{lem}		
}
%
\subsection{Proof or Proposition \ref{P:variance_y_0_non_zero}}
Throughout the proof we suppress in the notation the dependence upon $T$ of $h_1(T)$ and $h_2(T)$. The constant $C$ may change from line to line and is independent of $T$.
{\modar In the proof, we will use repeatedly Lemmas \ref{L:HShort}--\ref{L:HGlobal}. To this end, we consider a compact set $K$ that contains a ball of radius $\sqrt{2}$ centered at $(x_0,y_0)$, and as a result the support of $(x,y) \mapsto \varphi(x-x_0,y-y_0)$ is included in this compact.}

To prove the proposition, it is sufficient to prove that the following inequalities holds both, {\modar for $T$ large enough},
\begin{align} \label{E:maj_var_1}
&\Var( \hat{\pi}_T(x_0,y_0) )
\le C \frac{1}{T h_2} +  C \varepsilon(T,h_1(T),h_2(T)),
\\ \label{E:maj_var_2}
&\Var( \hat{\pi}_T(x_0,y_0) )
\le C \frac{1}{T \sqrt{h_1}} +  C \varepsilon(T,h_1(T),h_2(T)).
\end{align} 

\underline{First step:}  we prove \eqref{E:maj_var_1}.

From \eqref{E:Def_pi_hat}, and the stationarity of the process we get that 
\begin{equation} \label{E:Def_kappa}
\Var \left( \hat{\pi}_T(x_0,y_0) \right) = \frac{1}{T^2} \int_0^T \int_0^T \kappa(t-s) \dd t \dd s,
\end{equation}
where 
\begin{equation*}
\kappa(u)= \Cov \left( \varphi_{h_1,h_2}(X_0-x_0,Y_0-y_0) ,  \varphi_{h_1,h_2}(X_u-x_0,Y_u-y_0)\right).
\end{equation*}
We deduce that
\begin{equation} \label{E:control_Var_int}
\Var \left( \hat{\pi}_T(x_0,y_0) \right) \leq \frac{1}{T} \int_0^T \abs{\kappa(s)} \dd s.
\end{equation}
We will find an upper bound for the integral on the right-hand side of the latter expression by splitting the time interval $[0,T]$ into 4 pieces $[0,T]=[0,\delta) \cup [\delta,D_1) \cup [D_1,D_2) \cup {\modar [D_2,T]}$, where $\delta$, $D_1$, $D_2$, will be chosen latter.

$\bullet$ For $s \in [0,\delta)$, we write from \eqref{E:Def_kappa} and using Cauchy-Schwarz inequality and the stationarity of the process, 
\begin{align*}
\abs{\kappa(s)} &\le \Var ( \varphi_{h_1,h_2}(X_0-x_0,Y_0-y_0))^{1/2} \Var ( \varphi_{h_1,h_2}(X_s-x_0,Y_s-y_0))^{1/2} 
\\
&=\Var ( \varphi_{h_1,h_2}(X_0-x_0,Y_0-y_0)).
\end{align*}
This variance is smaller than 
\begin{equation*}
\int_{\mathbb{R}^2} \varphi_{h_1,h_2}(x-x_0,y-y_0)^2 \pi(x,y) \dd x \dd y 
\end{equation*}
and using that $\pi$ is {\modar bounded} and \eqref{E:def_var_phi_double}, we deduce
\begin{equation} \label{E:kappa_crude_control}
\abs{\kappa_s} \le  \frac{C}{h_1 h_2}.
\end{equation}
 In turn, we have
\begin{equation}
\label{E:int_kappa_morceau1}
\int_0^\delta \abs{\kappa(s)} \dd s \le C\frac{\delta}{h_1 h_2}.
\end{equation}

$\bullet$ For $s \in [\delta,D_1)$, where  $\delta < D_1  \leq 1$.
We write 
\begin{multline*}
\abs{\kappa(s)}\le
E\left[ \abs{\varphi_{h_1,h_2}(X_0-x_0,Y_0-y_0)} \abs{\varphi_{h_1,h_2}(X_u-x_0,Y_u-y_0)} \right]+ \\
E\left[ \abs{\varphi_{h_1,h_2}(X_0-x_0,Y_0-y_0) }\right]
E\left[ \abs{\varphi_{h_1,h_2}(X_s-x_0,Y_s-y_0) }\right]
\end{multline*}
Using that $(X_t,Y_t)_t$ is stationary with marginal law having a {\modar bounded} density we deduce that $E\left[ \abs{\varphi_{h_1,h_2}(X_s-x_0,Y_s-y_0) }\right] \leq C
\int_{\mathbb{R}^2} \abs{\varphi_{h_1,h_2}(x-x_0,y-y_0)} \dd x \dd y \le C$ from \eqref{E:def_var_phi_double}.
This gives,
\begin{multline} \label{E:maj_kappa_as_int}
\abs{\kappa(s)}\le 
\int_{\mathbb{R}^2} \abs{\varphi_{h_1,h_2}(x-x_0,y-y_0) }
\\
\int_{\mathbb{R}^2} \abs{\varphi_{h_1,h_2}(x'-x_0,y'-y_0) } 
p_s(x,y;x',y') \dd x' \dd y' \pi(x,y) \dd x \dd y
+ C.
\end{multline}
Using {\modar now equation \eqref{E:Hshort} in Lemma \ref{L:HShort},} we get that 
\begin{equation}
\label{E:decoupe_kappa_morceau2}
\abs{\kappa(s)} \le \overline{\kappa}^1(s) + \overline{\kappa}^2(s) + C
\end{equation}
 with
\begin{multline} \label{E:def_overline_kappa1}
\overline{\kappa}^1(s):=
\int_{\mathbb{R}^2} \abs{\varphi_{h_1,h_2}(x-x_0,y-y_0) }
\\
\int_{\mathbb{R}^2} \abs{\varphi_{h_1,h_2}(x'-x_0,y'-y_0) } 
p^{G}_s(x,y;x',y') \dd x' \dd y' \pi(x,y) \dd x \dd y,
\end{multline}  
\begin{multline*}
\overline{\kappa}^2(s):=
\int_{\mathbb{R}^2} \abs{\varphi_{h_1,h_2}(x-x_0,y-y_0) }
\\
\int_{\mathbb{R}^2} \abs{\varphi_{h_1,h_2}(x'-x_0,y'-y_0) } 
p^{U}_s(x,y;x,y) \dd x' \dd y' \pi(x,y) \dd x \dd y.
\end{multline*}  
In order to upper bound $\overline{\kappa}^1(s)$ we show that the Gaussian kernel \eqref{E:Def_p_G} appearing in the expression of $\overline{\kappa}^1(s)$ takes small values for $s \in [\delta,D_1)$ as soon as $\delta$ is well chosen. 
Recall that $y_0 \neq 0$, and for simplicity assume that $y_0 >0$. Then, using that $\varphi$ is compactly supported {\modar on$ [-1,1]$} we know that  $\abs{\varphi_{h_1,h_2}(x-x_0,y-y_0) \varphi_{h_1,h_2}(x'-x_0,y'-y_0)} \neq 0$ implies that 
\begin{equation}\label{E:def_rectangle}
\abs{x-x_0} \le h_1,\abs{x'-x_0} \le h_1 , \abs{y-y_0} \le h_2, \abs{y'-y_0} \le  h_2.
\end{equation}
 Let us denote {\modar $K(h_1,h_2)$} the rectangle of $\mathbb{R}^4$ defined by the conditions \eqref{E:def_rectangle}. Then,
\begin{multline} \label{E:cont_kappa1}
\overline{\kappa}^1(s) \leq C \int_{K(h_1,h_2)}
\abs{\varphi_{h_1,h_2}(x-x_0,y-y_0)}
	\\
	\abs{\varphi_{h_1,h_2}(x'-x_0,y'-y_0)} 
		p^{G}_s(x,y;x',y') \dd x \dd y \dd x' \dd y',
\end{multline} 
where we used that $\pi$ is {\modar bounded.}
On $K(h_1,h_2)$, we have $\frac{y+y'}{2} \ge \frac{y_0}{2}>0$ if $h_2$ is small enough, and $\abs{x-x'} \le 2 h_1$. Hence, if we assume that $s \ge  {\modar \frac{ 6 h_1}{y_0}} $ we have
$\abs{x'-x} \le \frac{s y_0}{3}$. It entails, $x'-x-\frac{y+y'}{2}s \le \frac{s y_0}{3} - \frac{y_0}{2}s= {\modar -\frac{sy_0}{6}}$, 
 and in turn $\frac{(x'-x-\frac{y+y'}{2}s)^2}{s^3} \ge {\modar  \frac{ y_0^2}{36} \frac{1}{s}}$. 
Plugging in \eqref{E:Def_p_G} this yields to $p^{G}_s(x,y;x';y') \le \frac{C}{s^2} \exp(-\frac{1}{C s})$ for some constant $C$ independent of $s$.  
Using \eqref{E:cont_kappa1}, we deduce
\begin{multline*}
\overline{\kappa}^1(s) \leq 
\frac{C}{s^2} \exp(-\frac{1}{C s})
\int_{K(h_1,h_2)}
\abs{\varphi_{h_1,h_2}(x-x_0,y-y_0)}
\\
\abs{\varphi_{h_1,h_2}(x'-x_0,y'-y_0)} 
\dd x \dd y \dd x' \dd y'
\end{multline*}
and hence 
\begin{equation} \label{E:cont_kappa1_bis}
\overline{\kappa}^1(s) \le 
\frac{C}{s^2} \exp(-\frac{1}{C s}).
\end{equation}

To control $\overline{\kappa}^2(s)$, we use \eqref{E:pte_p_U} and 
$\norme{
\varphi_{h_1,h_2}(\cdot-x_0,\cdot-y_0)}_\infty \le \frac{C}{h_1 h_2}$ to get that for all $x,y$ in {\modar the compact $K$}
containing a ball of radius $\sqrt{2}$ centered at $(x_0,y_0)$,  we have the upper bound
$\int_{\mathbb{R}^2}  \abs{\varphi_{h_1,h_2}(x'-x_0,y'-y_0) } 
p^{U}_s(x,y;x',y') \dd x' \dd y' \le \frac{C}{h_1h_2}e^{-\frac{1}{Cs}}$. As a consequence, 
\begin{align} \nonumber
\overline{\kappa}^2(s)
&\le C\int_{\mathbb{R}^2} \abs{\varphi_{h_1,h_2}(x-x_0,y-y_0)}
	\frac{1}{h_1 h_2} e^{-\frac{1}{Cs}} \pi(x,y) \dd x \dd y
	\\ \label{E:cont_kappa2}
& \le 	\frac{C}{h_1 h_2} e^{-\frac{1}{Cs}},
\end{align}
where we used again that $\pi$ is {\modar bounded and that the support of $\varphi_{h_1,h_2}(\cdot-x_0,\cdot-y_0)$ is included in $K$}.
From \eqref{E:decoupe_kappa_morceau2}, \eqref{E:cont_kappa1_bis}--\eqref{E:cont_kappa2}, we deduce that for $ {\modar \frac{6 h_1}{y_0} } \le \delta \le D_1 \le 1$,
\begin{align} \nonumber
\int_\delta^{D_1} \abs{\kappa(s)} \dd s
&\le \int_\delta^{D_1} [ \frac{C}{s^2} \exp(-\frac{1}{C s})+\frac{C}{h_1 h_2} \exp(-\frac{1}{C s}) + C] \dd s
\\ \nonumber
& \le \int_\delta^{D_1} [ \frac{C}{s^2} \exp(-\frac{1}{C s})+\frac{C }{h_1 h_2 s^2} \exp(-\frac{1}{C s}) + C] \dd s
\\ \label{E:int_kappa_morceau2}
& \le  C \exp(-\frac{1}{C D_1}) [1+ \frac{1}{h_1 h_2}] + C D_1 ,
\end{align}
where in the second line we used $D_1 \le 1$.

$\bullet$ For $s \in [D_1,D_2)$ with $D_1 \le 1 \le D_2 < T$, we start from the control \eqref{E:maj_kappa_as_int} 
{\modar that we write
\begin{equation*}
\abs{\kappa(s)}\le 
\int_{\mathbb{R}^2} \abs{\varphi_{h_1,h_2}(x-x_0,y-y_0) } P_s(\abs{\varphi_{h_1,h_2}(\cdot-x_0,\cdot-y_0)} )(x,y) \pi(x,y) \dd x \dd y
+ C.
\end{equation*}
}
{\modar Since 
$\varphi_{h_1,h_2}(\cdot-x_0,\cdot-y_0)$ vanishes outside the compact neighbourhood $K$ of $(x_0,y_0)$ we can use 
Lemma \ref{L:HGlobal} to upper bound the semi group term.}
Hence, we get for some constant $C>0$, 
{\modar 
\begin{multline*}
\abs{\kappa(s)} \le
C 
\int_{\mathbb{R}^2} \abs{\varphi_{h_1,h_2}(x-x_0,y-y_0) } 
\times
\\
[\frac{\norme{\varphi_{h_1,h_2}}_{L^1(\mathbb{R}^2)}}{D_1^2} + \frac{1}{h_1 h_2} e^{-1/(C D_1)}]
 \pi(x,y) \dd x \dd y 
+ C,
\end{multline*}  
}
and we deduce {\modar $\abs{\kappa(s)} \le C [\frac{1}{D_1^2} + \frac{1}{h_1 h_2} e^{-1/(C D_1)} +1 ] $}.
This yields, 
\begin{equation} \label{E:int_kappa_morceau3}
\int_{D_1}^{D_2} \abs{\kappa(s)} \dd s \le
 {\modar C [\frac{D_2}{D_1^2} + \frac{D_2}{h_1 h_2} e^{-1/(C D_1)} + D_2],}
\end{equation}
for some constant $C>0$.

$\bullet$ For $s \in [D_2,T]$, we use the covariance control \eqref{E:Mixing}, that allows us to write
$\abs{\kappa(s)} \le C \norme{\varphi_{h_1,h_2}(\cdot-x_0,\cdot-y_0)}_\infty^2 e^{-\rho s} \le 
C \left(\frac{1}{h_1h_2}\right)^2 e^{-\rho s}$, for $C>0$ and $\rho>0$. It entails the upper bound,
\begin{equation} \label{E:int_kappa_morceau4}
\int_{D_2}^T \abs{\kappa(s)} \dd s \le C \frac{e^{-\rho D_2}}{(h_1h_2)^2}.
\end{equation}

Collecting together \eqref{E:control_Var_int}, \eqref{E:int_kappa_morceau1}, \eqref{E:int_kappa_morceau2}, \eqref{E:int_kappa_morceau3}, \eqref{E:int_kappa_morceau4} we deduce,
\begin{multline*}
\Var(\hat{\pi}_T(x_0,y_0)) \le \frac{C}{T}\big[ \frac{\delta}{h_1 h_2} + \exp(-\frac{1}{C D_1}) [1+{\modar \frac{D_2}{h_1 h_2}}] + 
\\D_1 + {\modar D_2 + \frac{D_2}{D_1^2}}  + \frac{e^{-\rho D_2}}{(h_1h_2)^2}\big],
\end{multline*}
for $C>0$ some constant.
We choose $\delta= {\modar \frac{6 h_1}{y_0}}$, $D_1=\frac{1}{C \abs{\log h_1 h_2}}$, 
$D_2=\frac{\abs{\ln( (h_1 h_2)^2)}}{\rho}$. By \eqref{E:h_sub_poly_1}--\eqref{E:h_sub_poly_2} we see that this choice is such that,
$\frac{6 h_1 }{y_0}=\delta < D_1 < {\modar 1 <} D_2 <T$ for $T$ large enough. And it yields,
\begin{equation*}
\Var(\hat{\pi}_T(x_0,y_0)) \le
\frac{C}{T}\big[ \frac{{\modar 6 }}{y_0} \frac{1}{h_2}+ h_1h_2 + 
 D_1 + {\modar  D_2 + \frac{D_2}{D_1^2} } +1 \big].
\end{equation*}
Since $h_1 \to 0$, $h_2 \to 0$, and as the value of $C$ may change from line to line, we can write that
\begin{equation*}
\Var(\hat{\pi}_T(x_0,y_0)) \le \frac{C}{T} \left[ \frac{1}{h_2}+ \abs{\ln(h_1 h_2)}^C \right]
\end{equation*}
and we have shown \eqref{E:maj_var_1}.

\underline{Second step:} we prove \eqref{E:maj_var_2}.

We use the same decomposition of $\int_0^T \abs{\kappa(s)} \dd s$ in four terms as for the proof of \eqref{E:maj_var_1}, but we treat in a different way the contribution of the short time correlations $\int_0^\delta \abs{\kappa(s)} \dd s$. 

$\bullet$ Let us find an upper bound for $\kappa(s)$ for $s \in (0,\delta)$ with $\delta <1$.
We recall that, from {\modar Lemma \ref{L:HShort},}
the decomposition \eqref{E:decoupe_kappa_morceau2} holds true where $\overline{\kappa}^1(s)$ is given by \eqref{E:def_overline_kappa1} and $\overline{\kappa}^2(s)$ is  upper bounded by \eqref{E:cont_kappa2}. We now study $\overline{\kappa}^1(s)$. To this end, we 
remark that $p_s^G(x,y;x',y') \le \frac{C}{\sqrt{s}} q_s(x'\mid x,y,y')$ where
\begin{equation*}
 q_s(x'\mid x,y,y') = \frac{C}{s^{3/2}} \exp \left( -C^{-1} \frac{(x'-x-\frac{y+y'}{2} s )^2}{s^3}  \right).
\end{equation*}
Let us stress that 
\begin{equation} \label{E:prop_variance_temp_q}
\sup_{s \in (0,1)} \sup_{(x,y,y')\in \mathbb{R}^3}
\int_{\mathbb{R}} q_s(x'\mid x,y,y') \dd x' \le C < \infty.
\end{equation}
Thus, using \eqref{E:def_overline_kappa1}, we have
\begin{multline*}
\overline{\kappa}^1(s)
\le \frac{C}{\sqrt{s}}
\int_{\mathbb{R}^2} \abs{\varphi_{h_1,h_2}(x-x_0,y-y_0) } \pi(x,y) \\
\Big(
\int_{\mathbb{R}^2} \abs{\varphi_{h_1,h_2}(x'-x_0,y'-y_0) } 
 q_s(x'\mid x,y,y') \dd x' \dd y'\Big) \dd x \dd y.
\end{multline*} 
By \eqref{E:def_var_phi_double}, we have $\abs{\varphi_{h_1,h_2}(x'-x_0,y'-y_0) } \le \frac{C}{h_1} \frac{1}{h_2}\abs{\varphi(\frac{y'-y_0}{h_2})}$, and thus, using \eqref{E:prop_variance_temp_q}, we get
\begin{multline*}
\int_{\mathbb{R}^2} \abs{\varphi_{h_1,h_2}(x'-x_0,y'-y_0) } 
q_s(x'\mid x,y,y') \dd x' \dd y' 
\\ 
 \le \frac{C}{h_1} \int_{\mathbb{R}} \frac{1}{h_2}\abs{\varphi(\frac{y'-y_0}{h_2})}
\int_{\mathbb{R}}q_s(x'\mid x,y,y') \dd x'
  \dd y'
\\
 \le \frac{C}{h_1} \int_{\mathbb{R}} \frac{1}{h_2}\abs{\varphi(\frac{y'-y_0}{h_2})} \dd y'
\le  \frac{C}{h_1}.
\end{multline*}
We deduce that $\overline{\kappa}^1(s) \le \frac{C}{\sqrt{s} h_1}  
\int_{\mathbb{R}^2} \abs{\varphi_{h_1,h_2}(x-x_0,y-y_0) } \pi(x,y) \dd x \dd y \le \frac{C}{\sqrt{s} h_1}  $.
Collecting the latter equation with \eqref{E:decoupe_kappa_morceau2} and \eqref{E:cont_kappa2}, it yields 
\begin{align} \nonumber
\int_0^\delta \abs{\kappa(s)} \dd s & \le \int_0^\delta [\frac{C}{\sqrt{s}h_1} + C \frac{ e^{-\frac{1}{Cs}} }{h_1 h_2} + C] \dd s  
\\ \nonumber
& \le \int_0^\delta [\frac{C}{\sqrt{s}h_1} + C \frac{1}{s^2}\frac{ e^{-\frac{1}{Cs}} }{h_1 h_2} + C] \dd s
\\ \label{E:int_kappa_morceau1_alternative}
& \le C [ \frac{\sqrt{\delta}}{h_1} +   \frac{ e^{-\frac{1}{C\delta}} }{h_1 h_2} +\delta  ] 
\end{align}
where we have used in the second line that $\delta<1$.

We now gather \eqref{E:control_Var_int},  \eqref{E:int_kappa_morceau2}, \eqref{E:int_kappa_morceau3}, \eqref{E:int_kappa_morceau4}, \eqref{E:int_kappa_morceau1_alternative}, to derive,
\begin{multline*}
\Var(\hat{\pi}_T(x_0,y_0)) \le \frac{C}{T}\bigg[ \frac{\sqrt{\delta}}{ h_1} + \frac{ e^{-\frac{1}{C\delta}} }{h_1 h_2}+\delta+ 
\exp(-\frac{1}{C D_1}) [1+{\modar \frac{D_2}{h_1 h_2}}] + 
\\D_1 + {\modar  D_2+ \frac{D_1}{D_2^2}} + \frac{e^{-\rho D_2}}{(h_1h_2)^2}\bigg].
\end{multline*}
We choose the same thresholds as in the first step,  $\delta= \frac{6 h_1}{y_0}$, $D_1=\frac{1}{C \abs{\log h_1 h_2}}$, 
$D_2=\frac{\abs{\ln( (h_1 h_2)^2)}}{\rho}$. Recalling \eqref{E:h_sub_poly_1}--\eqref{E:h_sub_poly_2}, we have $e^{-(C 6 h_1 /y_0)^{-1}}=O(e^{- \varepsilon {\modar \log(T)^3}})=o(h_1h_2)$, with some $\varepsilon>0$.
We derive that
\begin{equation*}
\Var(\hat{\pi}_T(x_0,y_0)) \le \frac{C}{T}[\frac{1}{\sqrt{h_1}} + \abs{\ln(h_1h_2)}^C],
\end{equation*}
for some $C>0$. The second step of the proposition is proved.
\qed
\begin{rem}
\begin{enumerate}
	\item The Proposition \ref{P:variance_y_0_non_zero} consists actually in the two upper bounds \eqref{E:maj_var_1}--\eqref{E:maj_var_2} for the variance of the estimator. We see that one of these two bounds is smaller than the other, depending on the relative positions of $h_2$ or $\sqrt{h_1}$. It explains why the expression for the rate of convergence of the estimator in Theorem \ref{T:upper_bound} depends on the relative positions of $k_1$ and $k_2/2$, which determines which one of the two bounds \eqref{E:maj_var_1} or \eqref{E:maj_var_2} is used in the bias/variance decomposition of the estimation error (see proof of Theorem \ref{T:upper_bound}).
	\item The control of $\kappa(s)= \Cov \left( \varphi_{h_1,h_2}(X_0-x_0,Y_0-y_0) ,  \varphi_{h_1,h_2}(X_s-x_0,Y_s-y_0)\right)$ for $s \in [\delta,D_1]$, with $\delta \approx h_1$ and $D_1 \le 1$ depends on the fine structure of the main term \eqref{E:Def_p_G} in the short time expansion of the transition density of the process and on the fact that $y_0 \neq 0$. In the situation $y_0 \neq 0$, it is impossible to get such a refined result on the covariance, and eventually the bound on the variance of the estimator is larger (see Proposition \ref{P:variance_y_0_zero}). 
\end{enumerate}
\end{rem}
\subsection{Proof of Proposition \ref{P:variance_y_0_zero}}
We need to prove that the following two inequalities hold true, {\modar for $T$ large enough:}
\begin{align}
 \label{E:maj_var_1_y0_0}
&\Var( \hat{\pi}_T(x_0,y_0) )
\le C \frac{1}{T h_1^{2/3}} +  C \varepsilon(T,h_1(T),h_2(T)),
\\ \label{E:maj_var_2_y0_0}
&\Var( \hat{\pi}_T(x_0,y_0) )
\le C \frac{\ln (T)}{T h_2^2} +  C \varepsilon(T,h_1(T),h_2(T)).
\end{align}
{\modar Again we consider $K$ a compact set of $\mathbb{R}^2$ that contains a ball of radius $\sqrt{2}$ centered at $(x_0,y_0)$.}

\underline{First step :} let us prove \eqref{E:maj_var_1_y0_0}.

We recall the control \eqref{E:control_Var_int} for the variance of $\hat{\pi}_T(x_0,y_0)$ and split the integral in \eqref{E:control_Var_int} into four pieces corresponding to the partition
$[0,T]=[0,\delta) \cup [\delta,D_1) \cup [D_1,D_2) \cup {\modar [D_2,T]}$, where $\delta$, $D_1$, $D_2$ will be specified latter. 
Let us stress that in the proof of Proposition \ref{P:variance_y_0_non_zero}, only the control of $\abs{\kappa(s)}$ for $s\in [\delta,D_1)$, uses the fact that $y_0 \neq 0$.

$\bullet$ For $s \in [0,\delta)$ with $\delta<1$, we recall the result obtained in \eqref{E:int_kappa_morceau1_alternative} which states
\begin{equation*}
\int_0^\delta \abs{\kappa(s)} \dd s \le C 
[ \frac{\sqrt{\delta}}{h_1} +   \frac{ e^{-\frac{1}{C\delta}} }{h_1 h_2} +\delta  ].
\end{equation*}

$\bullet$ For $s \in [\delta,D_1]$ with $0<\delta<D_1<1$, exactly with the same proof as in Proposition \ref{P:variance_y_0_non_zero}, we have $\abs{\kappa(s)} \le \overline{\kappa}^1(s)+\overline{\kappa}^2(s)+C$
where $\overline{\kappa}^1(s)$ is given by \eqref{E:def_overline_kappa1} and $\overline{\kappa}^2(s)$ is upper-bounded as in \eqref{E:cont_kappa2}. We need to find a control on $\overline{\kappa}^1(s)$ in the situation $y_0 \neq 0$. Using, from \eqref{E:Def_p_G}, that $p^G_s(x,y;x',y) \le \frac{C}{s^2}$, {\modar $\forall (x,y,x',y') \in K^2$}, and the fact that $\pi$ is {\modar bounded,} we get
\begin{multline*}
\overline{\kappa}^1(s) \leq 
\frac{C}{s^2} 
\int_{\mathbb{R}^4}
\abs{\varphi_{h_1,h_2}(x-x_0,y-y_0)}
\\
\abs{\varphi_{h_1,h_2}(x'-x_0,y'-y_0)} 
\dd x \dd y \dd x' \dd y' \le \frac{C}{s^2}.
\end{multline*}
We deduce
\begin{align}\nonumber
\int_{\delta}^{D_1} \abs{\kappa(s)} \dd s &\le \int_{\delta}^{D_1}  C[ \frac{1}{s^2} + \frac{e^{-\frac{1}{Cs}}}{h_1h_2} + C] \dd s
\\
& \le \nonumber
C[\frac{1}{\delta}+ \exp(-\frac{1}{C D_1}) \frac{1}{h_1h_2}+D_1]
\\  \label{E:int_kappa_morceau2_y_0_0}
& \le
 C[\frac{1}{\delta}+ \exp(-\frac{1}{C D_1}) \frac{1}{h_1h_2}], 
\end{align}
where we used $D_1 \le 1$.

$\bullet$ For $s \in [D_1,D_2)$, with $D_1<1<D_2<T$, we use the control \eqref{E:int_kappa_morceau3}.	

$\bullet$ For $s \in [D_2,T]$, with $D_2<T$, we use the control \eqref{E:int_kappa_morceau4}.	
	
Collecting the four previous controls, we deduce
\begin{multline*}
\Var( \hat{\pi}_T(x_0,y_0))\le
\frac{C}{T} \big[ 
\frac{\sqrt{\delta}}{h_1}+ \frac{e^{-\frac{1}{C\delta}}}{h_1h_2}+ \frac{1}{\delta} + \frac{\exp(-\frac{1}{C D_1}) }{h_1h_2}{\modar(1+D_2)}+
		\\
		{\modar D_1+D_2+\frac{D_2}{D_1^2}} +\frac{e^{-\rho D_2}}{(h_1h_2)^2}
\big]
\end{multline*}	
where $C>0$, $\rho>0$. We choose $\delta$ that balances $\sqrt{\delta}/h_1$ with $1/\delta$, namely $\delta=h_1^{2/3}$ which is smaller than $1$ for $T$ large enough, recalling \eqref{E:h_sub_poly_2}. Next, we choose $D_2=\frac{\abs{\ln((h_1h_2)^2)}}{\rho}$, and $D_1=\frac{C}{\abs{\log{h_1 h_2}}}$. We deduce
\begin{equation*}
\Var( \hat{\pi}_T(x_0,y_0))\le \frac{C}{T}[ \frac{1}{h_1^{2/3}} + \abs{\log(h_1h_2)}^C  ],
\end{equation*}
for some $C>0$, and where we have used that, from \eqref{E:h_sub_poly_1}--\eqref{E:h_sub_poly_2},
 $\exp (-1/(Ch_1^{2/3}))= O( \exp( - \varepsilon {\modar \log(T)^{2}}))=o(h_1h_2)$, with some $\varepsilon >0$.	
 Hence \eqref{E:maj_var_1_y0_0} is proved.
 
\underline{Second step:} we show \eqref{E:maj_var_2_y0_0}.

Comparing with the first part of the proposition, we have to modify our upper bound on $\int_0^\delta \abs{\kappa(s)} \dd s$ with  $\delta<1$. We split this integral into two parts,
$\int_0^{\delta'} \abs{\kappa(s)} \dd s + \int_{\delta'}^{\delta} \abs{\kappa(s)} \dd s$. On the first part, we use the control \eqref{E:kappa_crude_control} and get
\begin{equation} \label{E:int_kappa_morceau_1_cut_left}
\int_0^{\delta'} \abs{\kappa(s)} \dd s \le C\frac{\delta'}{h_1 h_2}.
\end{equation}

On the second part, we use \eqref{E:decoupe_kappa_morceau2}, where $\overline{\kappa}^2(s)$ is bounded in \eqref{E:cont_kappa2}. We deduce,
\begin{align} \nonumber
\int_{\delta'}^{\delta} \abs{\kappa(s)} \dd s &\le \int_{\delta'}^{\delta}  C [\frac{1}{h_1 h_2} e^{-\frac{1}{C s}} + 1] \dd s
+ \int_{\delta'}^{\delta} \abs{\overline{\kappa}^1(s)} \dd s
\\ \label{E:int_kappa_morceau_1_cut_inter}
& \le C [\frac{1}{h_1 h_2} e^{-\frac{1}{C \delta}} + \delta]
+ \int_{\delta'}^{\delta} \overline{\kappa}^1(s) \dd s .
\end{align} 
To upper bound $\overline{\kappa}^1(s)$, we use \eqref{E:def_overline_kappa1} and \eqref{E:def_var_phi_double}, and obtain by Fubini's Theorem,
\begin{multline} \label{E:bound_kappa_1_y0_0_altenate}
\overline{\kappa}^1(s)=
	\int_{\mathbb{R}^2} \abs{\frac{1}{h_1}\varphi(\frac{x-x_0}{h_1}) \frac{1}{h_1}\varphi(\frac{x'-x_0}{h_1})}   
	\\	\big( \int_{\mathbb{R}^2} \abs{ \frac{1}{h_2}\varphi(\frac{y-y_0}{h_2}) \frac{1}{h_2}\varphi(\frac{y'-y_0}{h_2})} p^{G}_s(x,y;x',y') \pi(x,y)   \dd y \dd y' \big) \dd x \dd x'.
\end{multline}
Since {\modar $\pi$ is bounded,} and using \eqref{E:Def_p_G}, we deduce that the inner integral is lower than
\begin{align*} 
&\frac{\norme{\varphi}_\infty^2}{h_2^2}
\int_{\mathbb{R}^2}
\frac{C}{s^2} \exp \left( 
- \frac{1}{C} \left[ \frac{(y-y')^2}{s} + \frac{(x'-x-\frac{y+y'}{2}s)^2}{s^3} \right]
\right) \dd y \dd y'
\\
=&
\frac{\norme{\varphi}_\infty^2}{2h_2^2}
\int_{\mathbb{R}^2}
\frac{C}{s} \exp \left( 
- \frac{1}{C} \left[ w^2 + (\frac{x'-x}{s^{3/2}}-\frac{w'}{2})^2 \right]
\right) \dd w \dd w'
\\
=&
\frac{\norme{\varphi}_\infty^2}{2h_2^2}
\int_{\mathbb{R}^2}
\frac{C}{s} \exp \left( 
- \frac{1}{C} \left[ w^2 + (\frac{w'}{2})^2 \right]
\right) \dd w \dd w'  ,
\end{align*}	
where we have made the change of variables $w=\frac{y-y'}{\sqrt{s}}$, $w'=\frac{y+y'}{\sqrt{s}}$ in the second line, and used the invariance by translation of the Lebesgue measure in the last one.
We deduce that the inner integral in \eqref{E:bound_kappa_1_y0_0_altenate} is lower than
$\frac{C}{h_2^2 s}$ where we stress that $C$ does not depend on $(x,x')$. In turn, 
\begin{equation*}
\overline{\kappa}^1(s)\le \frac{C}{h_2^2 s}
\int_{\mathbb{R}^2} \abs{\frac{1}{h_1}\varphi(\frac{x-x_0}{h_1}) \frac{1}{h_1}\varphi(\frac{x'-x_0}{h_1})} \dd x \dd x' \le \frac{C}{h_2^2 s}.
\end{equation*}
 This yields, using
\eqref{E:int_kappa_morceau_1_cut_inter} to 
\begin{equation}
\label{E:int_kappa_morceau_1_cut_right}
\int_{\delta'}^{\delta} \abs{\kappa(s)} \dd s \le C [\frac{1}{h_1 h_2} e^{-\frac{1}{C \delta}} + \delta+ \frac{1}{h_2^2} \ln(\frac{\delta}{\delta'} )].
\end{equation}
Collecting \eqref{E:int_kappa_morceau_1_cut_left}, \eqref{E:int_kappa_morceau_1_cut_right}, \eqref{E:int_kappa_morceau2_y_0_0}, \eqref{E:int_kappa_morceau3}, \eqref{E:int_kappa_morceau4}, we deduce, for $0<\delta'<\delta<D_1<1<D_2<T$,
\begin{multline*}
\Var(\hat{\pi}_T(x_0,y_0)) \le \frac{C}{T}\big[ \frac{\delta'}{h_1h_2}+ 
 \frac{e^{-\frac{1}{C \delta}}}{h_1 h_2}  +  \frac{1}{h_2^2} \ln(\frac{\delta}{\delta'} ) + \frac{1}{\delta} + \frac{\exp(-\frac{1}{C D_1}) }{h_1h_2} {\modar (1+D_2)}+
\\
{\modar D_1+ D_2 + \frac{D_1}{D_2^2}} 
+\frac{e^{-\rho D_2}}{(h_1h_2)^2} \big].
\end{multline*}
Now we let $\delta=h_2^2$ and $\delta'=T^{-C_\delta'}$ where $C_\delta'>0$ is such that
$\delta'=o(h_1h_2)$ 
which is possible from the at most polynomial decay of the bandwidths, resorting to \eqref{E:h_sub_poly_1}.
As in the first step of the proposition, we set $D_1=\frac{C}{\abs{\ln(h_1h_2)}}$, $D_2=\frac{\abs{\ln((h_1h_2)^2)}}{\rho}$. With these choices, we have for $T$ large enough,
$0<\delta'<\delta<D_1<1<D_2<T$ and
\begin{align*}
\Var(\hat{\pi}_T(x_0,y_0)) 
&\le \frac{C}{T}[1+ \frac{e^{-\frac{1}{Ch_2^2}}}{h_1h_2}+\frac{1}{h_2^2} \ln(\delta/T^{-C_{\delta'}})+\abs{\ln(h_1h_2)}^{C}]
\\
&\le \frac{C}{T}[\frac{\ln(T)}{h_2^2} +\abs{\ln(h_1h_2)}^{C}]
\end{align*}
where we used again \eqref{E:h_sub_poly_1}--\eqref{E:h_sub_poly_2} in the last line. This proves
\eqref{E:maj_var_2_y0_0}.
\qed
{\modar 
\subsection{Numerical simulations}
In this section, we explore
numerically on an example the behaviour of the variance of the estimator as $h_1$ and $h_2$ go to $0$. Especially, we wonder
if the variance of the estimator asymptotically depends, in these simulations, only on 
  the minimum of two quantities related to $h_1$ and $h_2$ as suggested by the upper bounds in Propositions \ref{P:variance_y_0_non_zero}--\ref{P:variance_y_0_zero}. This is the crucial point in the upper bound of the variance,   that makes the choice of the optimal bandwidth very specific, allowing an arbitrary thin bandwidth on one component.
  We consider the model \eqref{E:model_X}--\eqref{E:model_Y}
   with $\beta(x,y)=0.5$, $V(x)=x^2/2$ and $a(x,y)=1$. From a Mont\'e--Carlo experiment based on 500 replications, we evaluate the variance of $\hat{\pi}(x_0,y_0)$ for $(x_0,y_0)=(0,1.5)$ and with different values of bandwidths $h_1$ and $h_2$. We have chosen the simple kernel $\varphi(u)=\frac{1}{2} 1_{[-1,1]}(u)$ and $T=200$. 
  Results are given in Figure \ref{F:variance}, where each curve corresponds to a choice for the bandwidth $h_2$, and these curves plot the value of the variance as a function of $h_1$, using log-scales. We see that, as expected, the variance is increasing as $h_1$ and $h_2$ get smaller.
  Moreover, it appears that when $h_2$ gets smaller than some threshold depending on $h_1$ the variance ceases to strictly increase, as all the curves are flat on the left side of the Figure \ref{F:variance}.  Symmetrically,
  we see that the right part of the curves for $h_2=10^{-1.8}$ and $h_2=10^{-2.4}$ coincides.
   It shows for instance that decreasing $h_2$ below $10^{-1.8}$ does not increase anymore the variance when $h_1 \ge 10^{-2}$. Hence, the numerical results are consistent with the upper bound given for the variance, 
  as a function depending on   $\min (1/h_2,1/\sqrt{h_1})$. This suggests that the upper bounds given in Propositions \ref{P:variance_y_0_non_zero}--\ref{P:variance_y_0_zero} are fairly sharp.
  
  \begin{figure}[h]\caption{Variance of the estimator}
  	\label{F:variance}
  \includegraphics[width=\textwidth]{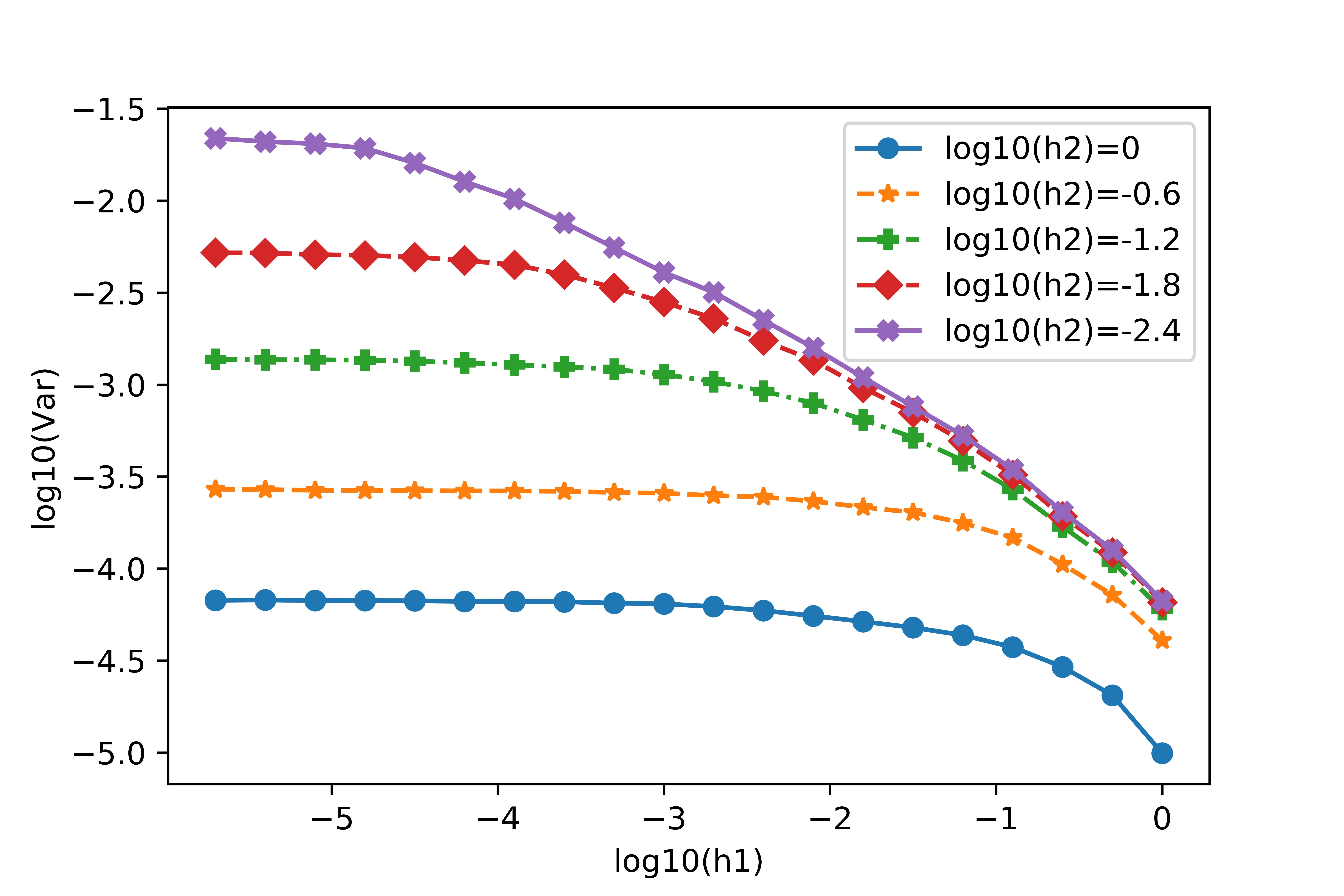}
  \end{figure}
   
%
%

}
 
\section{Minimax lower bound} \label{S:minimax}
In this section, we show that it is impossible to construct any estimator with a {\modar uniform}  rate better (up to a log term) 
 than the rates obtained in
Theorems \ref{T:upper_bound}--\ref{T:upper_bound_y0_0}.

\subsection{Lower bounds}
For the computation of lower bounds, we introduce the family of S.D.E. 
\begin{align}
\label{E:inf_bound_mod_X}
&\dd X_t = Y_t \dd t \\
\label{E:inf_bound_mod_Y}
&\dd Y_t = 2 \sigma  \dd B_t - [ \sigma^2 \beta(X_t,Y_t) Y_t + V'(X_t)] \dd t,
\end{align}
where $\sigma>0$, $\beta$ is {\modar a bounded $\mathcal{C}^1$ function} lower bounded by a strictly positive number and $V$ is {\modar $\mathcal{C}^2$.}
As the model satisfies the conditions of \HReg{} we know that the S.D.E. admits a weak solution.
We know that if $V$ satisfies \HErg{} then a Lyapounov function exists and the process admits a  unique stationary measure, that we note $\pi_{V,\beta}$. In Section \ref{Ss:link_drift_statio}, we make more explicit the connection between $\pi_{V,\beta}$ and the coefficients $V$, $\beta$.
 Remark that we omit in the notations the dependence on $\sigma$, as $\sigma$ will be fixed in the sequel.

If the stationary measure exists and is unique,  we denote $\mathbb{P}_{V,\beta}$ the law of a stationary solution $(X_t,Y_t)_{t \ge 0} $ of \eqref{E:inf_bound_mod_X}--\eqref{E:inf_bound_mod_Y}. Here, $\mathbb{P}_{V,\beta}$ is a measure on the space of continuous function $\mathcal{C}([0,\infty),\mathbb{R}^2)$, and we note by $\mathbb{E}_{V,\beta}$ the corresponding expectation. When needed, we will note
by $\mathbb{P}_{V,\beta}^{(T)}$ the law of the stationary process $(X_t,Y_t)_{t \in [0,T]} $ solution to the S.D.E. \eqref{E:inf_bound_mod_X}--\eqref{E:inf_bound_mod_Y}.

In the sequel, we note again $(X,Y)$ the canonical process on $\mathcal{C}([0,\infty),\mathbb{R}^2)$ or $\mathcal{C}([0,T],\mathbb{R}^2)$.

In order to write down an expression for the minimax risk of estimation, we have to consider a set of solutions to the S.D.E. \eqref{E:inf_bound_mod_X}--\eqref{E:inf_bound_mod_Y}, which are stationary and whose stationary measure has a prescribed H\"{o}lder regularity. This leads us to the following definition.

\begin{defi}\label{D:class_drift}
	 Let $V : \mathbb{R} \to \mathbb{R}$ be a  {\modar $\mathcal{C}^2$} function satisfying \HErg{}. We consider  $k_1>0,~k_2>0,~R>0$, and 
	 {\modar $R'>1$} real numbers. We define $\Sigma^{k_1,k_2}(V,R,R')$ the set of {\modar continuously differentiable} functions $\beta : \mathbb{R}^2 \to \mathbb{R}$ satisfying the following two conditions :
	 \begin{enumerate}
		\item {\modar $1/R' \le \beta(x,y) \le R'$} for all $(x,y) \in \mathbb{R}^2$,
		\item the density $\pi_{V,\beta}$ of the stationary measure associated to the S.D.E. \eqref{E:inf_bound_mod_X}--\eqref{E:inf_bound_mod_Y} 
		is such that
		 $\pi_{V,\beta} \in \mathcal{H}^{k_1,k_2}(R)$.
	\end{enumerate}
\end{defi}
We introduce the minimax risk for the estimation at some point. Let $(x_0,y_0) \in \mathbb{R}^2$, and $V$,  $k_1,~k_2,~R,~R'$ as in Definition \ref{D:class_drift}.
We let
\begin{equation} \label{E:minimax_risk}
R_T(V,k_1,k_2,R,R')
= \inf_{\widetilde{\pi}_T} \sup_{\beta \in \Sigma^{k_1,k_2}(V,R,R')}
\mathbb{E}_{V,\beta} \left[ (\widetilde{\pi}_T(x_0,y_0) - \pi_{V,\beta}(x_0,y_0))^2\right],
\end{equation}
where the infimum is taken on all possible estimators of $\pi_T(x_0,y_0)$, that is for $\widetilde{\pi}_T=\widetilde{\pi}_T(x_0,y_0)$ ranging in the set of all the measurable functions of $(X_t,Y_t)_{t \in [0,T]}$  with values in $\mathbb{R}$.

\begin{theo} \label{T:lower_bound_y0_non_zero}
Let $k_1,k_2,R>0$ and assume $y_0 \neq 0$ and $\max(k_1,k_2/2)>1/2$. Then, there exists $V$ satisfying \HErg{} and {\modar $R'>1$} such that, for some constant $C>0$, we have :
\begin{equation}\label{E:lower_bound_y0_non_zero}
R_T(V,k_1,k_2,R,R') \ge C T^{-2v(k_1,k_2)}, \quad \forall T>1,
\end{equation}	
with 
\begin{equation*}
v(k_1,k_2)=\begin{cases}
\frac{k_2}{2k_2+1}, &\text{ if $k_1 < k_2/2$,}
\\
\frac{k_1}{2k_1+1/2}, &\text{ if $k_1\ge k_2/2$.}
\end{cases}
\end{equation*}
\end{theo}
\begin{rem}
	\begin{enumerate}
		\item 
Theorem \ref{T:lower_bound_y0_non_zero} tells us that it is impossible to find an estimator with a rate of estimation, for the pointwise $L^2$ risk, better than $T^{-v(k_1,k_2)}$ {\modar on a the class of diffusions $Z=(X,Y)$ having a with $\mathcal{H}^{k_1,k_2}(R)$ stationary measure.}
On the other hand the estimator introduced in Section \ref{S:Estimator} achieves this rate, by Theorem \ref{T:upper_bound}, {\modar for each diffusion $Z=(X,Y)$ satisfying \HReg{} and \HErg{} and with stationary measure in $\mathcal{H}^{k_1,k_2}(R)$.}

\item 
	{\modar  
		The upper bound given in Theorem \ref{T:upper_bound} is not stated uniformly on the class of all diffusions satisfying \HReg{} and \HErg{} and is not a minimax upper bound. To get uniform upper bound, we would need that the mixing control \eqref{E:Mixing} holds uniformly on a class of diffusions whose coefficients satisfy uniform versions of the assumptions \HReg, \HErg. We are not aware of such uniform mixing results, and hence, getting a uniform version of Theorem \ref{T:upper_bound} is left for further research.}
 
	\item The condition $\max(k_1,k_2/2)>1/2$ asserts that $\pi$ is not too irregular with respect to both variables $x$ and $y$.
	Such assumption is weak, as the stationary measure is typically 
	smoother than the coefficients of the S.D.E. (see point 3 of Remark \ref{R:link_pi_beta} below).
	\end{enumerate}
\end{rem}

\begin{theo} \label{T:lower_bound_y0_zero}
	Let $k_1,k_2,R>0$ and assume $y_0 =0$ and $\max(k_1,k_2/3)>2/3$. Then, there exists $V$ satisfying \HErg{} and {\modar $R'>1$} such that, for some constant $C>0$, we have :
	\begin{equation} \label{E:lower_bound_y0_zero}
	R_T(V,k_1,k_2,R,R') \ge C T^{-2v'(k_1,k_2)}, \quad \forall T>1,
	\end{equation}	
	with 
	\begin{equation*}
	v'(k_1,k_2)=\begin{cases}
	\frac{k_2}{2k_2+2}, &\text{ if $k_1 < k_2/3$,}
	\\
	\frac{k_1}{2k_1+2/3}, &\text{ if $k_1\ge k_2/3$.} 
	\end{cases}
	\end{equation*}
\end{theo}	
Again, the previous result shows that the estimator introduced in Section \ref{S:Estimator} is rate efficient, up to a log term, in the case where $y_0=0$.

\subsection{Explicit link between the drift and the stationary measure} \label{Ss:link_drift_statio}
{\modar
Recall that from Proposition \ref{P:HErg_HMixv2}, \HReg{} and \HErg{} are sufficient for the existence and uniqueness of a stationary probability of the process solution of \eqref{E:inf_bound_mod_X}--\eqref{E:inf_bound_mod_Y} (see \cite{Wu01} also, or see Talay \cite{Talay} for related conditions too).
In this section, we will  characterize explicit relations between $(V,\beta)$ and $\pi_{V,\beta}$. 
}

 We need to introduce $A^\star_{V,\beta}$ the adjoint on $\mathbf{L}^2(\mathbb{R}^2,\dd x\dd y)$ of the generator $A_{V,\beta}$ of the process $Z=(X,Y)$ solution to \eqref{E:inf_bound_mod_X}--\eqref{E:inf_bound_mod_Y}. 

Assume that $V$ is {\modar $\mathcal{C}^2$} and that $(x,y) \mapsto y\beta(x,y) $ is of class {\modar $\mathcal{C}^{1}$.} Then we define for $g : \mathbb{R}^2 \to \mathbb{R}$ any {\modar $\mathcal{C}^{2}$} function,
\begin{multline} \label{E:def_A_star}
A^\star_{V,\beta}g(x,y)=2 \sigma^2 \frac{\partial^2 g}{\partial y^2}(x,y) - y \frac{\partial g}{\partial x}(x,y) 
\\
+
 [\sigma^2 y \beta(x,y)+V'(x)] \frac{\partial g}{\partial y}(x,y)
+
\sigma^2 \frac{\partial (y\beta)}{\partial y}(x,y) g(x,y).
\end{multline}
It can be checked that \eqref{E:def_A_star} is the expression for the adjoint of the generator of the process. If $g : \mathbb{R}^2 \to \mathbb{R}$ is a probability density, of class {\modar $\mathcal{C}^{2}$,}
solution to $A^\star_{V,\beta} g=0$, then it is an invariant density for the process. Hence, when the stationary distribution $\pi_{V,\beta}$ is unique, it can be computed as solution of the equation $A^\star_{V,\beta} \pi_{V,\beta}=0$. From the expression \eqref{E:def_A_star} it seems impossible to find explicit solutions $g$ to the equation $A^\star_{V,\beta} g=0$ for any $V$ and $\beta$, as one need to solve explicitly some P.D.E. 
Consequently, it seems impossible to write $\pi_{V,\beta}$ as an explicit expression of $(V,\beta)$.

On the other hand, using \eqref{E:def_A_star} it can be seen that if one consider $g$ and $V$ as fixed and $\beta$ as the unknown variable in the equation $A^\star_{V,\beta}g=0$, then finding solution in $\beta$ is simpler as one has to deal with a P.D.E. involving only differentiation with respect to $y$. As a consequence, it will be possible to express $\beta$ as a function of the stationary distribution $\pi$ (for a fixed $V$). This is the object of the next proposition. We need to introduce some notations first.

For $g \in \mathcal{C}^{2}$ and $g>0$, we define for all $(x,y) \in \mathbb{R}^2$,
\begin{equation}\label{E:def_xi}
\xi_g(x,y)=\frac{1}{\sigma^2 g(x,y)} \int_0^y 
\big[ z \frac{\partial g}{\partial x}(x,z) -V'(x) \frac{\partial g}{\partial y}(x,z)+
 2 \sigma^2 \frac{\partial^2 g}{\partial y^2} (x,z) 
\big] \dd z,
\end{equation}
and 
\begin{equation} \label{E:def_beta_g}
\beta_g(x,y)=
\begin{cases} 
\displaystyle
\frac{1}{y} \xi_g (x,y), \text{ for $y\neq 0$,} 
\\ \displaystyle
\lim_{y\to 0}\frac{1}{y} \xi_g(x,y)=
\\ \quad  \displaystyle
\frac{1}{\sigma^2 g(x,0)}[-V'(x) \frac{\partial g}{\partial y}(x,0)+
2 \sigma^2 \frac{\partial^2 g}{\partial y^2} (x,0)],  \text{ for $y = 0$.}
\end{cases}
\end{equation}

\begin{prop}\label{P:link_pi_beta}
\begin{enumerate}
	\item[1)]
Let $V : \mathbb{R} \to \mathbb{R}$ with regularity {\modar $\mathcal{C}^2$,} $g : \mathbb{R}^2 \to \mathbb{R}$ with regularity 
{\modar $\mathcal{C}^{2}$} 
and $g >0$.

Then,
we have that $(x,y) \mapsto y \beta_g(x,y)$ is a {\modar $\mathcal{C}^{1}$} function and
\begin{equation} \label{E:A_star_beta_g_zero}
A^\star_{V,\beta_g} g(x,y)=0, \quad \forall (x,y) \in \mathbb{R}^2.
\end{equation}
Moreover, $\beta_g$ is the unique function solution to \eqref{E:A_star_beta_g_zero} such that $(x,y) \mapsto y \beta_g(x,y)$ is 
{\modar  $\mathcal{C}^{0,1}$.}

\item[2)] Let $V : \mathbb{R} \to \mathbb{R}$ with regularity {\modar $\mathcal{C}^2$} and satisfying \HErg{} and consider $\pi : \mathbb{R}^2 \to \mathbb{R}$ a probability density with regularity {\modar $ \mathcal{C}^{2}$} and $\pi >0$. 

Assume that {\modar $1/R'<\beta_\pi < R'$ for some $R'>1$,} where $\beta_\pi$ is defined by \eqref{E:def_beta_g}.

Then, $\pi$ is the unique stationary probability of the S.D.E \eqref{E:inf_bound_mod_X}--\eqref{E:inf_bound_mod_Y} with damping coefficient $\beta=\beta_\pi$ and potential $V$
\end{enumerate}	
\end{prop}
\begin{pf}
1) For $g : \mathbb{R}^2 \to \mathbb{R}$ with regularity {\modar $ \mathcal{C}^{2}$} and $g >0$ and $\beta$ such that $(x,y) \mapsto y\beta(x,y)$ is of class {\modar $\mathcal{C}^{1}$,} we can write the equation $A^\star_{V,\beta}g(x,y)=0$, recalling \eqref{E:def_A_star}, as
\begin{equation} \label{E:ODE_beta}
 \frac{\partial (y \beta) }{\partial y}(x,y) + \frac{\frac{\partial g}{\partial y}(x,y) }{g(x,y)} (y \beta(x,y))=
i_g(x,y)  
\end{equation}
with 
\begin{equation*}
i_g(x,y)=\frac{1}{\sigma^2 g(x,y)} \left[
y \frac{\partial g }{\partial x} (x,y) -V'(x) \frac{\partial g }{\partial y} (x,y) - 2 \sigma^2 \frac{\partial^2 g}{\partial y^2}(x,y)
\right].
\end{equation*}
Let us fix $x \in \mathbb{R}$, we then interpret \eqref{E:ODE_beta} as an ordinary differential equation with differentiation variable $y$, where the unknown parameter is the function
$y \mapsto \xi(x,y):=y\beta(x,y)$ :
\begin{equation}
\label{E:ODE_xi}
 \frac{\partial \xi }{\partial y} + \frac{\frac{\partial g}{\partial y} }{g} \xi =
i_g.  
\end{equation}
 A solution of the homogeneous equation $\frac{\partial \xi }{\partial y} + \frac{\frac{\partial g}{\partial y} }{g} \xi =0$ is $\xi(x,y)=1/g(x,y)$. Then, by variation of the constant method,
we deduce that the solution to \eqref{E:ODE_xi}, 
has the expression,
\begin{equation*}
\xi(x,y)=\frac{1}{g(x,y)}
\left[
c(x)+\int_0^y g(x,z) i_g(x,z) \dd z
\right],
\end{equation*}
where $c(x)$ is an integration constant. As $\xi(x,0)=0 \times \beta(x,0)=0$, we deduce that $c(x)=0$, $\forall x$. Hence the solution $\xi(x,y)=y \beta(x,y)$ of \eqref{E:ODE_xi} 
is given by \eqref{E:def_xi} and in turn, we deduce that $\beta=\beta_g$ given by \eqref{E:def_beta_g} is the unique solution to \eqref{E:ODE_beta} or equivalently to \eqref{E:A_star_beta_g_zero}.

2) Using Ito's formula, one can check that any $\pi$ solution to $A^\star_{V,\beta} \pi=0$ is a stationary measure for the process $(X,Y)$ given by \eqref{E:inf_bound_mod_X}--\eqref{E:inf_bound_mod_Y}. 
From the first part of the proposition, $\pi$ is solution to $A^\star_{V,\beta_\pi} \pi=0$. 
By {\modar Proposition \ref{P:HErg_HMixv2},}
the stationary measure of the equation with damping coefficient $\beta_\pi$ is unique, and is thus equal to $\pi$. 
\end{pf}
\begin{rem} \label{R:link_pi_beta}
	\begin{enumerate}
	\item Is is known that the S.D.E. 
	\begin{align*}
	&\dd X_t = Y_t \dd t \\
	&\dd Y_t = 2 \sigma  \dd B_t - [ \sigma^2 \varepsilon Y_t + V'(X_t)] \dd t.
	\end{align*}
	admits for the stationary measure $\pi(x,y) = C \exp \left( -\frac{\varepsilon}{2} [
	\frac{ y^2}{2} + V(x) ] \right)$ (see e.g. \cite{Comte_et_al_17}).  
	As expected, if we take $\pi(x,y)= C \exp \left( -\frac{\varepsilon}{2} [
	\frac{y^2}{2} + V(x) ] \right)$ and compute $\beta_\pi$ by the formula \eqref{E:def_beta_g}, we find $\beta_\pi=\varepsilon$.
	\item Proposition \ref{P:link_pi_beta} shows how to compute the damping part of the drift in order to get a diffusion with a prescribed stationary measure. However, it is not clear that for a given $\pi$ the corresponding $\beta_\pi$, computed with \eqref{E:def_xi}--\eqref{E:def_beta_g} satisfies the sign condition $\beta_\pi>1/R'$ insuring that the process is indeed ergodic. This is why in part 2) of Proposition \ref{P:link_pi_beta} we postulate  $\beta_\pi >1/R'$. However, we will see in Section \ref{S:proof_lower_bound} that if $\pi$ is a small deviation of $\pi_0$ given by   $\pi_0(x,y)= C \exp \left( -\frac{\varepsilon }{2} [
	\frac{y^2}{2} + V(x) ] \right)$, then the corresponding $\beta$ is a small deviation of $\beta_0=\varepsilon$ and thus is positive.
	\item The equation \eqref{E:def_xi} enables to relate the degree of smoothness of the drift coefficient and of the stationary measure, when the latter exists and is unique. 
	Indeed, from   \eqref{E:def_xi}, we get that if $\pi \in \mathcal{C}^{k_1,k_2}$ and $V \in \mathcal{C}^{k_1}$, then the associated drift of the S.D.E. $\xi_{\pi} + V'$ is $\mathcal{C}^{k_1-1,k_2-1}$. 
	\end{enumerate}
\end{rem}
\subsection{Proof of Theorem \ref{T:lower_bound_y0_non_zero}} \label{S:proof_lower_bound}
The proof of the lower bound is made by a comparison between the minmax risk \eqref{E:minimax_risk} and  some Bayesian risk where the Bayesian prior is supported on a set of two elements.
\subsubsection{Construction of the prior} \label{Ss:Constru_prior_non_zero}
Let $k_1,k_2>0$ and $R>0$. We set $V_0(x)=\abs{x}^2$ and define
\begin{equation}\label{E:def_pi_0_beta_0}
\pi_0(x,y)=c_\eta \exp( -\frac{\eta }{2} [ \frac{y^2}{2}+ x^2  ]), \quad  \beta_0(x,y)=\eta, \quad  \xi_0(x,y)=\eta y,
\end{equation}
where $\eta>0$ and
where $c_\eta$ is the constant that makes $\pi_0$ a probability measure. 
The function $\pi_0$ is $\mathcal{C}^\infty$ and it is possible to choose {\modar $0<\eta<1/2$} small enough such that
\begin{equation*}
\pi_0 \in \mathcal{H}^{k_1,k_2}(R/2).
\end{equation*}
We know from Section \ref{Ss:link_drift_statio} that $\pi_0$ is the unique stationary measure for $(X^{(0)},Y^{(0)})$ solution of
\begin{align}
\label{E:def_X_0}
&\dd X^{(0)}_t = Y^{(0)}_t \dd t \\
\label{E:def_Y_0}
&\dd Y^{(0)}_t = 2 \sigma \dd B_t - [ \sigma^2 
\eta Y^{(0)}_t+ V_0'(X^{(0)}_t)] \dd t.
\end{align}
Now, if we set ${\modar R'=2/\eta>1}$, we have, using $\beta_0 = \eta$ and recalling Definition \ref{D:class_drift},
\begin{equation} \label{E:beta_0_in_set}
\beta_0 \in \Sigma^{k_1,k_2}(V,R/2,R'/2).
\end{equation} 
Let $h : \mathbb{R} \to \mathbb{R}$ be a $\mathcal{C}^\infty$ function with support on $[-1,1]$ and such that 
\begin{equation} \label{E:pty_h}
h(0)=1, ~ \int_{-1}^1 h(z) \dd z=0, ~ \int_{-1}^1 zh(z) \dd z=0.
\end{equation}
 We set for $T>0$,
\begin{equation}\label{E:def_tilde_pi_prior}
\tilde{\pi}_T(x,y)= \pi_0(x,y) + \frac{1}{M_T} h(\frac{x-x_0}{h_1(T)}) h(\frac{y-y_0}{h_2(T)}),
\end{equation}
where  $M_T$, $h_1(T)$, $h_2(T)$ will be calibrated later and satisfy
\begin{equation*}
M_T \xrightarrow{T\to \infty} \infty,~ 
h_1(T) \xrightarrow{T\to \infty} 0,~ 
h_2(T) \xrightarrow{T\to \infty} 0.
\end{equation*}
From \eqref{E:pty_h}, we see that $\int_{\mathbb{R}^2} \tilde{\pi}_T(x,y) \dd x \dd y= \int_{\mathbb{R}^2} \pi_0(x,y) \dd x \dd y=1$, and using $\pi_0>0$, $1/M_T \to 0$ and that $h$ is compactly supported, we see that $\tilde{\pi}_T >0$ for $T$ large enough. Hence $\tilde{\pi}_T$ is a smooth probability measure for $T$ large enough. We define 
\begin{equation*}
\tilde{\beta}_T(x,y)= \beta_{\tilde{\pi}_T}(x,y), \quad  \tilde{\xi}_T(x,y) = y \tilde{\beta}_T(x,y)=\xi_{\tilde{\pi}_T}(x,y),
\end{equation*}
where we used the definitions \eqref{E:def_xi} and \eqref{E:def_beta_g}.

Before proving Theorem \ref{T:lower_bound_y0_non_zero}, we need to state two lemmas. 
The first lemma shows that the two functions $\beta_0$ and $\tilde{\beta}_T$ only differ on some vanishing neighbourhood of $(x_0,y_0)$.

\begin{lem}\label{L:beta_beta_T}
1) Let us define the compact set of $\mathbb{R}^2$
\begin{equation*}
K_T=[x_0-h_1(T),x_0+h_1(T)]\times [y_0-h_2(T),y_0+h_2(T)].
\end{equation*}
Then, for $T$ large enough, we have for all $(x,y) \notin K_T$ :
\begin{equation*}
\beta_0(x,y)=\tilde{\beta}_T(x,y),\quad \xi_0(x,y)=\tilde{\xi}_T(x,y).
\end{equation*}
2) For $(x,y) \in K_T$, we have the control
\begin{align}
\label{E:maj_diff_beta}
&\abs{ \beta_0(x,y)-\tilde{\beta}_T(x,y)}  \le \frac{C}{M_T}\left[\frac{h_2(T)}{h_1(T)}+ \frac{1}{h_2(T)}\right],
\\ \label{E:maj_diff_xi}
&\abs{\xi_0(x,y)-\tilde{\xi}_T(x,y)}  \le \frac{C}{M_T}\left[\frac{h_2(T)}{h_1(T)}+ \frac{1}{h_2(T)}\right],
\end{align}
where $C$ is some constant independent of $T$, $h_1(T)$, $h_2(T)$, $M_T$.

3) We have
\begin{equation*}
\int_{\mathbb{R}^2} \abs{ \tilde{\xi}_T(x,y) - \xi_0(x,y)     }^2 \dd x \dd y
\le 
\frac{C}{M_T^2} \left[\frac{h_2(T)^3}{h_1(T)}+ \frac{h_1(T)}{h_2(T)}\right].
\end{equation*}
\end{lem}
\begin{pf}
1) We first prove the $\tilde{\xi}_T$ and $\xi_0$ coincides on $K_T^c$. With the definition \eqref{E:def_xi} in mind, we set for $g$ of class $\mathcal{C}^{1,2}$:
\begin{align*}
\mathcal{I}[g](x,y)&=
\frac{1}{\sigma^2} \int_0^y 
\big[ z \frac{\partial g}{\partial x}(x,z) -V'(x) \frac{\partial g}{\partial y}(x,z)+
2 \sigma^2 \frac{\partial^2 g}{\partial y^2} (x,z) 
\big] \dd z
\\
&=\mathcal{I}^1[g](x,y)+\mathcal{I}^2[g](x,y)+\mathcal{I}^3[g](x,y),
\end{align*}
where
\begin{align}
\label{E:def_I_1_g}
\mathcal{I}^1[g](x,y)&= \frac{1}{\sigma^2} \int_0^y 
 z \frac{\partial g}{\partial x}(x,z) \dd z,
\\
\label{E:def_I_2_g}
\mathcal{I}^2[g](x,y)&= - \frac{V'(x)}{\sigma^2} [g(x,y)-g(x,0)],
\\
\label{E:def_I_3_g}
\mathcal{I}^3[g](x,y)&=2 [\frac{\partial g}{\partial y}(x,y) -\frac{\partial g}{\partial y}(x,0) ]. 
\end{align}
Using this notation, we have
\begin{equation} \label{E:rel_xi_I}
\tilde{\xi}_T=\frac{1}{\tilde{\pi}_T} \mathcal{I}[\tilde{\pi}_T], \quad 
\xi_0=\frac{1}{\pi_0} \mathcal{I}[\pi_0].
\end{equation}
Let us note 
\begin{equation}
\label{E:diff_pi}
d_T=\tilde{\pi}_T-\pi_0
\end{equation}
 and by \eqref{E:def_tilde_pi_prior}, we have
\begin{equation} \label{E:def_d_T}
d_T(x,y)=\frac{1}{M_T} h(\frac{x-x_0}{h_1(T)})h(\frac{y-y_0}{h_2(T)}).
\end{equation}
Since $g  \mapsto \mathcal{I}[g]$ is  a linear operator we deduce that
\begin{equation} \label{E:rel_xi_tilde_I_linear}
\tilde{\xi}_T=\frac{1}{\tilde{\pi}_T} \mathcal{I}[\tilde{\pi}_T]=\frac{1}{\tilde{\pi}_T} \mathcal{I}[\pi_0] +
\frac{1}{\tilde{\pi}_T} \mathcal{I}[d_T].
\end{equation}
If $(x,y) \notin K_T$ we have from \eqref{E:diff_pi}, \eqref{E:def_d_T} and the fact that the support of $h$ is included in $[-1,1]$ that $\tilde{\pi}_T(x,y)=\pi_0(x,y)$. Thus $\tilde{\xi}_T(x,y)=\frac{1}{\pi_0(x,y)} \mathcal{I}[\pi_0](x,y) +
\frac{1}{\pi_0(x,y)} \mathcal{I}[d_T](x,y)=\xi_0(x,y)+
\frac{1}{\pi_0(x,y)} \mathcal{I}[d_T](x,y)$. It follows that the equality of $\tilde{\xi}_T$ and $\xi_0$ on $K_T^c$ will a consequence of the following fact:
\begin{equation} \label{E:I_d_T_outside}
\text{for $(x,y) \notin K_T$, we have,} \quad  \mathcal{I}[d_T](x,y)=0.
\end{equation}
Let us check that \eqref{E:I_d_T_outside} holds true. To this end, it is enough that $\mathcal{I}^{i}[d_T](x,y)=0$ for $i=1,2,3$ and $(x,y) \notin K_T$. Since $h$ is a smooth function with compact support on $[-1,1]$, the function $d_T$ and its derivatives vanishes outside of the compact set
$K_T$ by \eqref{E:def_d_T}. Recalling $y_0 \neq 0$, for $T$ large enough and for all $x \in \mathbb{R}$, the point $(x,0)$ does not belong to $K_T$, thus we deduce from \eqref{E:def_I_2_g}--\eqref{E:def_I_3_g} that $\mathcal{I}^{2}[d_T](x,y)=\mathcal{I}^{3}[d_T](x,y)=0$ when $(x,y) \notin K_T$. It remains to see that $\mathcal{I}^{1}[d_T](x,y)=0$ for $(x,y) \notin K_T$. We have
by \eqref{E:def_I_1_g} and \eqref{E:def_d_T},
\begin{equation}\label{E:I_1_d_T}
\mathcal{I}^{1}[d_T](x,y)=\frac{h'((x-x_0)/h_1(T))}{\sigma^2 M_T h_1(T)}\int_0^y z h (\frac{z-y_0}{h_2(T)}) \dd z.
\end{equation}
For $(x,y) \notin K_T$, a first possibility is $x \notin [x_0-h_1(T),x_0+h_1(T)]$ that leads to $\mathcal{I}^{1}[d_T](x,y)=0$ as $h'$ vanishes outside $[-1,1]$ and thus $h'((x-x_0)/h_1(T))=0$. 
Otherwise, we must have $y \notin [y_0-h_2(T),y_0+h_2(T)]$. For simplicity of the presentation, assume that $y_0>0$. Then,
\begin{align*}
\int_0^y z h (\frac{z-y_0}{h_2(T)}) \dd z
&=
\int_0^y z h (\frac{z-y_0}{h_2(T)}) 
1_{[y_0-h_2(T),y_0+h_2(T)]}(z) \dd z,
\\
&=
\begin{cases}
\displaystyle 0, &\text{  if $y \displaystyle \le y_0-h_2(T)$,}
 \\
\displaystyle \int_{y_0-h_2(T)}^{y_0+h_2(T)} z h (\frac{z-y_0}{h_2(T)}) \dd z, &\text{ if $ \displaystyle y \ge y_0+h_2(T)$.}
\end{cases}
\end{align*}
But  $\int_{y_0-h_2(T)}^{y_0+h_2(T)} z h (\frac{z-y_0}{h_2(T)}) \dd z = h_2(T)\int_{-1}^1 (y_0+h_2(T)z)h(z) \dd z=0$ using \eqref{E:pty_h}. This yields to
$\mathcal{I}^{1}[d_T](x,y)=0$ for $(x,y) \notin K_T$ and \eqref{E:I_d_T_outside} is proved. It follows that $\tilde{\xi}(x,y)=\xi_0(x,y)$ for $(x,y) \notin K_T$. 

The equality between $\tilde{\beta}_T$ and $\beta_0$ outside $K_T$ is a consequence of $\tilde{\beta}_T(x,y)=\tilde{\xi}(x,y)/y$, $\beta_0(x,y)={\xi_0}(x,y)/y$ for $y \neq 0$.

2) We will prove \eqref{E:maj_diff_xi} first. From \eqref{E:rel_xi_I} and \eqref{E:rel_xi_tilde_I_linear}, we have
\begin{equation}\label{E:diff_xi_tilde_xi_0}
\tilde{\xi}_T-\xi_0=\frac{\pi_0-\tilde{\pi}_T}{\tilde{\pi}_T} \xi_0+ \frac{1}{\tilde{\pi}_T}\mathcal{I}[d_T].
\end{equation}
On the set $K_T$, we see that $\tilde{\pi}_T=\pi_0 + d_T$ is lower bounded away from $0$ and that $\xi_0$ is bounded. Using $\norme{d_T}_{\infty} \le C/M_T$ we deduce that
\begin{equation*}
\forall (x,y) \in K_T,~ \abs{\tilde{\xi}_T(x,y)-\xi_0(x,y)} \le C \big[\frac{1}{M_T} + \mathcal{I}[d_T](x,y) \big].
\end{equation*}
Now, $\abs{\mathcal{I}[d_T](x,y)} \le \abs{\mathcal{I}^1[d_T](x,y)}+ \abs{ \mathcal{I}^2[d_T](x,y)} + \abs{ \mathcal{I}^3[d_T](x,y)}$. From \eqref{E:def_I_2_g} and \eqref{E:def_d_T},
$\abs{\mathcal{I}^2[d_T](x,y)} \le C \norme{d_T}_\infty \le C/M_T$ and by \eqref{E:def_I_3_g}, $\abs{\mathcal{I}^3[d_T](x,y)} \le C \norme{\frac{\partial d_T}{\partial y}}_\infty \le C/(h_2(T)M_T)$. 
Using \eqref{E:I_1_d_T}, we have 
\begin{equation*}
\abs{\mathcal{I}^1[d_T](x,y)} \le \frac{\norme{h'}_\infty}{\sigma^2 M_T h_1(T)} \int_{y_0-h_2(T)}^{y_0+h_2(T)}
\abs{z} \dd s \norme{h}_\infty \le C \frac{h_2(T)}{M_Th_1(T)}.
\end{equation*}
 We deduce that,
\begin{equation*}
\forall (x,y) \in K_T,~ \abs{\tilde{\xi}_T(x,y)-\xi_0(x,y)} \le \frac{C}{M_T} \big[ \frac{h_2(T)}{h_1(T)}+ 1 + \frac{1}{h_2(T)}    \big],
\end{equation*}
which gives \eqref{E:maj_diff_xi} as $h_2(T) \to 0$.

Eventually, \eqref{E:maj_diff_beta} follows from the fact that, for $T$ large enough, $K_T$ does not intersect the axis $y=0$ since $y_0 \neq 0$ and the relation between
 $\tilde{\beta}_T(x,y)=\tilde{\xi}(x,y)/y$, $\beta_0(x,y)={\xi_0}(x,y)/y$.

3) We have,
\begin{equation*}
\int_{\mathbb{R}^2} \abs{ \tilde{\xi}_T(x,y) - \xi_0(x,y)     }^2 \dd x \dd y
=
\int_{K_T} \abs{ \tilde{\xi}_T(x,y) - \xi_0(x,y)     }^2 \dd x \dd y
\end{equation*}
and the third point of the lemma follows from \eqref{E:maj_diff_xi}
with the fact that the Lebesgue measure of $K_T$ is proportional to $h_1(T)h_2(T)$. 
\end{pf}
\begin{lem} \label{L:beta_in_prior}
Let $\varepsilon >0$ and assume that for all $T$ large,
\begin{equation} \label{E:cond_reg_prior}
M_T^{-1} \le \varepsilon h_1(T)^{k_1}, \quad M_T^{-1} \le \varepsilon h_2(T)^{k_2},
\end{equation}
and
\begin{equation*}
\frac{h_2(T)}{h_1(T)} + \frac{1}{h_2(T)} = o(M_T) \text{ as $T \to \infty$.}
\end{equation*}
Then, if $\varepsilon>0$ is small enough, 
we have 
\begin{equation*}
\tilde{\beta}_T \in \Sigma^{k_1,k_2}(V_0,R,R'), 
\end{equation*} 
for all $T$ sufficiently large.	
\end{lem}
\begin{pf}
From Lemma \ref{L:beta_beta_T}, we know that $\tilde{\beta}_T=\beta_0$ outside $K_T$ and thus is constant equal to $\eta>0$ outside $K_T$. For $(x,y) \in K_T$, we have by \eqref{E:maj_diff_beta} in Lemma \ref{L:beta_beta_T},
{\modar  $\tilde{\beta}_T(x,y) = \beta_0(x,y) + O\left( \frac{1}{M_T}\left[\frac{h_2(T)}{h_1(T)}+ \frac{1}{h_2(T)} \right]\right) = \eta + o(1)$.}
	 where $C$ is some constant. 
	 Thus for $T$ sufficiently large we have
	 {\modar 
\begin{equation*}
\forall (x,y), \quad
1/R'=\eta/2<\tilde{\beta}_T(x,y) <1<R' 
\end{equation*}
where we recall that $R'=2/\eta>1$.
}
As $V_0$ is {\modar $\mathcal{C}^1$} and satisfy \HErg, we can apply the second point of Proposition \ref{P:link_pi_beta} and deduce that $\tilde{\pi}_T$ is the unique stationary measure associated to $\tilde{\beta}_T$. 
Recalling Definition \ref{D:class_drift}, the lemma will be shown as soon as we have,
\begin{equation*}
\tilde{\pi}_T \in \mathcal{H}^{k_1,k_2}(R).
\end{equation*}
Let us check the H\"older condition with respect to the variable $x$, as the condition with respect to the variable $y$ is similar. 
For all $(x,y) \in \mathbb{R}^2$ and $z \in [-1,1]$,
\begin{multline*}
\abs{\frac{\partial^{\entier{k_1}} \tilde{\pi}_T}{\partial x^{\entier{k_1}}}(x+z,y)-\frac{\partial^{\entier{k_1}} \tilde{\pi}_T}{\partial x^{\entier{k_1}}}(x,y)} 
\le 
\abs{\frac{\partial^{\entier{k_1}} \pi_0}{\partial x^{\entier{k_1}}}(x+z,y)-\frac{\partial^{\entier{k_1}} \pi_0}{\partial x^{\entier{k_1}}}(x,y)}+
\\
\abs{\frac{\partial^{\entier{k_1}} d_T}{\partial x^{\entier{k_1}}}(x+z,y)-\frac{\partial^{\entier{k_1}} d_T}{\partial x^{\entier{k_1}}}(x,y)}
\\
\le
\frac{R}{2} \abs{z}^{k_1-\entier{k_1}}+
\abs{\frac{\partial^{\entier{k_1}} d_T}{\partial x^{\entier{k_1}}}(x+z,y)-\frac{\partial^{\entier{k_1}} d_T}{\partial x^{\entier{k_1}}}(x,y)}
\\
\le
\frac{R}{2} \abs{z}^{k_1-\entier{k_1}}+ \frac{\norme{h}_\infty}{M_T h_1(T)^\entier{k_1}}
\abs{
h^{(\entier{k_1})}\left(\frac{x+z-x_0}{h_1(T)} \right)
- h^{(\entier{k_1})}\left(\frac{x-x_0}{h_1(T)} \right)
}
\end{multline*}
where we have successively used $\tilde{\pi}_T=\pi_0+d_T$, $\pi_0 \in \mathcal{H}^{k_1,k_2}(R/2)$, and the definition \eqref{E:def_d_T} of $d_T$.
We now write 
\begin{multline*}
\abs{
	h^{(\entier{k_1})}\left(\frac{x+z-x_0}{h_1(T)} \right)
	- h^{(\entier{k_1})}\left(\frac{x-x_0}{h_1(T)} \right)
}\le
\\
\abs{
h^{(\entier{k_1})}\left(\frac{x+z-x_0}{h_1(T)} \right)
- h^{(\entier{k_1})}\left(\frac{x-x_0}{h_1(T)} \right)
}^{k_1-\entier{k_1}}*(2\norme{h^{(\entier{k_1})}}_\infty)^{1-(k_1-\entier{k_1})}
\end{multline*}
which is smaller than
$\norme{h^{(\entier{k_1}+1)}}_\infty^{k_1-\entier{k_1}} \abs{\frac{z}{h_1(T)}}^{k_1-\entier{k_1}}* (2\norme{h^{(\entier{k_1})}}_\infty)^{1-(k_1-\entier{k_1})}$.
It implies that
\begin{equation*}
\abs{\frac{\partial^{\entier{k_1}} \tilde{\pi}_T}{\partial x^{\entier{k_1}}}(x+z,y)-\frac{\partial^{(\entier{k_1})} \tilde{\pi}_T}{\partial x^{\entier{k_1}}}(x,y)} \le \abs{z}^{k_1-\entier{k_1}}
\left[ \frac{R}{2} + \frac{c_h}{M_T h_1(T)^{k_1}}
\right]
\end{equation*}
where $c_h=\norme{h}_\infty \norme{h^{(\entier{k_1}+1)}}_\infty^{k_1-\entier{k_1}} (2\norme{h^{(\entier{k_1})}}_\infty)^{1-(k_1-\entier{k_1})}$. If one uses \eqref{E:cond_reg_prior} with any $\varepsilon < \frac{R}{2 c_h}$, we deduce
\begin{equation*}
\abs{\frac{\partial^{\entier{k_1}} \tilde{\pi}_T}{\partial x^{\entier{k_1}}}(x+z,y)-\frac{\partial^{\entier{k_1}} \tilde{\pi}_T}{\partial x^{\entier{k_1}}}(x,y)} \le R \abs{z}^{k_1-\entier{k_1}}.
\end{equation*}
This is the required H\"older control on the derivatives of $\tilde{\pi}_T$ with respect to $x$.	The lemma follows.	
\end{pf}
\subsubsection{Proof of the lower bound \eqref{E:lower_bound_y0_non_zero} on the minimax risk}
\label{Ss:proof_lower_bound_y_0_non_zero}
Let us recall some notations. We denote $\mathbb{P}_{V,\beta}$ the law of the stationary solution to \eqref{E:inf_bound_mod_X}--\eqref{E:inf_bound_mod_Y} on the canonical space $\mathcal{C}([0,\infty),\mathbb{R}^2)$ and $\mathbb{E}_{V,\beta}$ the corresponding expectation. We denote by $\mathbb{P}_{V,\beta}^{(T)}$ (resp. $\mathbb{E}_{V,\beta}^{(T)}$) the restrictions of this probability (resp. expectation) on $\mathcal{C}([0,T],\mathbb{R}^2)$.

Let $\widetilde{\pi}_T(x_0,y_0)$ be any measurable function from  $\mathcal{C}([0,T],\mathbb{R}^2)$ to $\mathbb{R}$. 
We will estimate by below, for $T$ large, 
\begin{equation*}
R(\tilde{\pi}_T(x_0,y_0)):=\sup_{\beta \in \Sigma^{k_1,k_2} (V_0,R,R')}
 \mathbb{E}^{(T)}_{V_0,\beta} 
 \left[ ( \widetilde{\pi}_T(x_0,y_0) - \pi_{V_0,\beta}(x_0,y_0) )^2  \right].
\end{equation*}

Let us assume that the following conditions hold true, 
\begin{align}
\label{E:cond_beta_set1_in_proof} 
&M_T^{-1} \le \varepsilon h_1(T)^{k_1}, \quad M_T^{-1} \le \varepsilon h_2(T)^{k_2},
\\
&\label{E:cond_beta_set2_in_proof}
\frac{h_2(T)}{h_1(T)} + \frac{1}{h_2(T)} = o(M_T) \text{ as $T \to \infty$,}
\end{align}
where $\varepsilon$ is sufficiently small to get the conclusion of Lemma \ref{L:beta_in_prior}. We deduce that for $T$ large enough $\tilde{\beta}_T \in \Sigma^{k_1,k_2}(V_0,R,R')$. From \eqref{E:beta_0_in_set}, we have $\beta_0 \in  \Sigma^{k_1,k_2}(V_0,R/2,R'/2) \subset \Sigma^{k_1,k_2}(V_0,R,R')$.
It follows
\begin{multline*}
R(\tilde{\pi}_T(x_0,y_0)) \ge \frac{1}{2} 
\mathbb{E}^{(T)}_{V_0,\tilde{\beta}_T} 
\left[ ( \widetilde{\pi}_T(x_0,y_0) - \pi_{V_0,\tilde{\beta}_T}(x_0,y_0) )^2  \right]
+ \\ \frac{1}{2}
\mathbb{E}^{(T)}_{V_0,\beta_0} 
\left[ ( \widetilde{\pi}_T(x_0,y_0) - \pi_{V_0,\beta_0}(x_0,y_0) )^2  \right].
\end{multline*}
Using Lemma \ref{L:Girsanov} below, we know that $Z^{(T)}=\frac{\dd \mathbb{P}_{V_0,\tilde{\beta}_T}^{(T)}}{\dd \mathbb{P}_{V_0,\beta_0}^{(T)}}$ exists, and we can write
\begin{align*}
R(\tilde{\pi}_T(x_0,y_0)) &\ge
 \frac{1}{2} 
\mathbb{E}^{(T)}_{V_0,\beta_0} 
\left[ ( \widetilde{\pi}_T(x_0,y_0) - \pi_{V_0,\tilde{\beta}_T}(x_0,y_0) )^2  Z^{(T)}\right]
+ \\  &
\quad \quad \quad \quad \quad \quad \quad 
\frac{1}{2}
\mathbb{E}^{(T)}_{V_0,\beta_0} 
\left[ ( \widetilde{\pi}_T(x_0,y_0) - \pi_{V_0,\beta_0}(x_0,y_0) )^2  \right]
\\
& 
\ge
\frac{1}{2 \lambda} 
\mathbb{E}^{(T)}_{V_0,\beta_0} 
\left[ ( \widetilde{\pi}_T(x_0,y_0) - \pi_{V_0,\tilde{\beta}_T}(x_0,y_0) )^2  1_{ \{ Z^{(T)} \ge \frac{1}{\lambda} \}} \right]
+ \\  &
\quad \quad \quad \quad \quad \quad \quad 
\frac{1}{2}
\mathbb{E}^{(T)}_{V_0,\beta_0} 
\left[ ( \widetilde{\pi}_T(x_0,y_0) - \pi_{V_0,\beta_0}(x_0,y_0) )^2  1_{ \{ Z^{(T)} \ge \frac{1}{\lambda} \}} \right]
\\
= 
\frac{1}{2 \lambda} 
\mathbb{E}^{(T)}_{V_0,\beta_0} 
&
\left[ [( \widetilde{\pi}_T(x_0,y_0) - \pi_{V_0,\tilde{\beta}_T}(x_0,y_0) )^2  + ( \widetilde{\pi}_T(x_0,y_0) - \pi_{V_0,\beta_0}(x_0,y_0) )^2 ]  1_{ \{ Z^{(T)} \ge \frac{1}{\lambda} \}} \right]
\end{align*}
for all $\lambda >1$.
 As
$( \widetilde{\pi}_T(x_0,y_0) - \pi_{V_0,\tilde{\beta}_T}(x_0,y_0) )^2  + 
( \widetilde{\pi}_T(x_0,y_0) - \pi_{V_0,\beta_0}(x_0,y_0) )^2  \ge \left( \frac{\pi_{V_0,\tilde{\beta}_T}(x_0,y_0)-\pi_{V_0,{\beta}_0}(x_0,y_0)}{2}  \right)^2 $ we deduce,
\begin{equation*}
R(\tilde{\pi}_T(x_0,y_0)) \ge
\frac{1}{8 \lambda} 
( \pi_{V_0,\tilde{\beta}_T}(x_0,y_0)-\pi_{V_0,{\beta}_0}(x_0,y_0) )^2 
 \mathbb{P}^T_{V_0,\beta_0}\left(Z^{(T)} \ge \frac{1}{\lambda}\right).
\end{equation*}
 Since $\pi_{V_0,\tilde{\beta}_T}=\tilde{\pi}_T$, $\pi_{V_0,{\beta}_0}=\pi_0$ and recalling \eqref{E:diff_pi}, \eqref{E:def_d_T} with $h(0)=1$ we deduce $\pi_{V_0,\tilde{\beta}_T}(x_0,y_0)-\pi_{V_0,{\beta}_0}(x_0,y_0)=d_T(x_0,y_0)=1/M_T$, and it follows,
\begin{equation*}
R(\tilde{\pi}_T(x_0,y_0)) \ge
\frac{1}{8\lambda} \frac{1}{M_T^2} \mathbb{P}^T_{V_0,\beta_0}\left(Z_T \ge \frac{1}{\lambda}\right).
\end{equation*}
From Lemma \ref{L:Girsanov} below
 we know that $\inf_{T \ge 0} \mathbb{P}^T_{V_0,\beta_0}\left(Z_T \ge \frac{1}{\lambda}\right)>0$ 
 for some $\lambda=\lambda_0$ 
 as soon as
\begin{equation*}
\sup_{T \ge 0} T \int_{\mathbb{R}^2}   \abs{\widetilde{\xi}_T(x,y)-\xi_0(x,y)}^2 \dd x \dd y    < \infty.
\end{equation*}
Using the third point of Lemma \ref{L:beta_beta_T}, the latter condition is implied by,
\begin{equation} \label{E:cond_Girsanov_in_proof}
\sup_T \frac{T}{M_T^2}  \left[\frac{h_2(T)^3}{h_1(T)}+ \frac{h_1(T)}{h_2(T)}\right] < \infty.
\end{equation}
We deduce that
\begin{equation} \label{E:lower_bound_generic}
R(\tilde{\pi}_T(x_0,y_0)) \ge
\frac{C}{M_T^2},
\end{equation}
for $C>0$, if the conditions \eqref{E:cond_beta_set1_in_proof}, \eqref{E:cond_beta_set2_in_proof} and \eqref{E:cond_Girsanov_in_proof} are satisfied. It remains to find the larger choice for $1/M_T^2$, subject to the conditions
\eqref{E:cond_beta_set1_in_proof}, \eqref{E:cond_beta_set2_in_proof} and \eqref{E:cond_Girsanov_in_proof}. The optimal choice depends on $k_1$ and $k_2$.

\underline{Case 1, $k_1<k_2/2$ :} 

We set $h_1(T)=h_2(T)^2$, and $h_2(T)=\left( \frac{1}{ \varepsilon M_T} \right)^{1/k_2}$. The choice for $h_2(T)$ saturates one the conditions in
\eqref{E:cond_beta_set1_in_proof}. Let us see that the other condition holds also true . Indeed $1/M_T=\varepsilon h_2(T)^{k_2}=
\varepsilon h_1(T)^{k_2/2}\le \varepsilon h_1(T)^{k_1}$ for $T$ large, as $k_1<k_2/2$ and $h_1(T) \to 0$. Thus \eqref{E:cond_beta_set1_in_proof} is satisfied.

Plugging the values of $h_1(T)$ and $h_2(T)$ in \eqref{E:cond_Girsanov_in_proof},
we obtain the constraint $\frac{T}{M_T^2} \left( \frac{1}{ \varepsilon M_T} \right)^{1/k_2} \le  C$ for some $C>0$, that leads us to the choice $M_T=T^{1/(2+1/k_2)}=T^{\frac{k_2}{2k_2+1}}$. Then, the condition \eqref{E:cond_beta_set2_in_proof} is satisfied as $k_2>1$, indeed
$h_2(T)/h_1(T)+1/h_2(T)=2/h_2(T)=2\left( \varepsilon M_T \right)^{1/k_2} =o(M_T)$.

Eventually, we deduce from the application of \eqref{E:lower_bound_generic},
\begin{equation} \label{E:res_lower_bound_case1_proof}
R(\tilde{\pi}_T(x_0,y_0)) \ge C/M_T^2=
C T^{-\frac{2k_2}{2k_2+1}}.
\end{equation}

\underline{Case 2, $k_1 \ge k_2/2$ :} 

We set $h_1(T)=h_2(T)^2$, and $h_1(T)=\left( \frac{1}{ \varepsilon M_T} \right)^{1/k_1}$. As $1/M_T=\varepsilon h_1(T)^{k_1} = h_2(T)^{2k_1} \le
h_2(T)^{k_2}$, we see that 
\eqref{E:cond_beta_set1_in_proof} is satisfied. Plugging these choices of bandwidths in \eqref{E:cond_Girsanov_in_proof}, we obtain the constraint $\frac{T}{M_T^2} \left( \frac{1}{ \varepsilon M_T} \right)^{1/(2k_1)} \le  C$ for some $C>0$, that leads us to the choice $M_T=T^{1/(2+1/(2k_1))}=T^{\frac{k_1}{2k_1+1/2}}$. Then, the condition \eqref{E:cond_beta_set2_in_proof} is satisfied as $k_1>1/2$, indeed
$h_2(T)/h_1(T)+1/h_2(T)=2/h_2(T)=2\left( \varepsilon M_T \right)^{1/(2k_1)} =o(M_T)$.

Eventually, we deduce from the application of \eqref{E:lower_bound_generic},
\begin{equation} \label{E:res_lower_bound_case2_proof}
R(\tilde{\pi}_T(x_0,y_0)) \ge C/M_T^2=
C T^{-\frac{2k_1}{2k_1+1/2}}.
\end{equation}
Gathering \eqref{E:res_lower_bound_case1_proof} and \eqref{E:res_lower_bound_case2_proof}, we have shown Theorem \ref{T:lower_bound_y0_non_zero}.  
\qed

\begin{lem} \label{L:Girsanov}
	1) The measure $\mathbb{P}^{(T)}_{V_0,\tilde{\beta}_T}$ is absolutely continuous with respect to $\mathbb{P}^{(T)}_{V_0,\beta_0}$.

2) Denote $Z^{(T)}=\frac{\dd \mathbb{P}^{(T)}_{V_0,\tilde{\beta}_T}}{\dd \mathbb{P}^{(T)}_{V_0,\beta_0}}$ and assume that,
\begin{equation} \label{E:cond_Girsanov_lemma_statement}
\sup_{T \ge 0} T \int_{\mathbb{R}^2}   \abs{\widetilde{\xi}_T(x,y)-\xi_0(x,y)}^2 \dd x \dd y    < \infty.
\end{equation}
Then, there exist $\lambda_0,~C>0$ such that,
\begin{equation*}
 \mathbb{P}^{(T)}_{V_0,\beta_0}
 \left(
 Z^{(T)} \ge \frac{1}{\lambda_0}
 \right) \ge C
\end{equation*}
for all $T$ large enough.
\end{lem}
\begin{pf}
1) The absolute continuity  $\mathbb{P}^{(T)}_{V_0,\tilde{\beta}_T} \ll \mathbb{P}^{(T)}_{V_0,\beta_0}$ and expression for the ratio $Z^{(T)}=\frac{\dd \mathbb{P}^{(T)}_{V_0,\tilde{\beta}_T}}{\dd \mathbb{P}^{(T)}_{V_0,\beta_0}}$ is obtained by Girsanov formula, changing the drift of the component $Y^0$ in \eqref{E:def_Y_0} to the drift appearing in the component $Y^{(T)}$ of the stationary solution of the S.D.E.
\begin{align*}
&\dd X^{(T)}_t = Y^{(T)}_t \dd t \\
&\dd Y^{(T)}_t = 2 \sigma \dd B_t - [ \sigma^2 \tilde{\beta}_T(Y^{(T)}_t) Y^{(T)}_t+ V_0'(X^{(T)}_t)] \dd t.
\end{align*}
By classical computations (see Theorem 1.12 in \cite{Kutoyants04}), we have
\begin{multline} \label{E:Girsanov_explicit}
Z^{(T)}=\frac{\dd \mathbb{P}^{(T)}_{V_0,\tilde{\beta}_T}}{\dd \mathbb{P}^{(T)}_{V_0,\beta_0}}
\left( (X_s,Y_s)_{0 \le s \le T}  \right)
\\
=\frac{\tilde{\pi}_T(X_0,Y_0)}{\pi_0(X_0,Y_0)}
\exp \Bigg\{
-\frac{1}{4}
\int_0^T \Delta_T(X_s,Y_s) \dd Y_s
\\
-\frac{1}{8 \sigma^2} \int_0^T [\sigma^2 {\xi}_0(X_s,Y_s) + \sigma^2 \Delta_T(X_s,Y_s) + V_0'(X_s,Y_s) ]^2 -  [\sigma^2 {\xi}_0(X_s,Y_s) + V_0'(X_s,Y_s) ]^2 \dd s
\Bigg\},
\end{multline}
where $\Delta_T(x,y)=\tilde{\xi}_T(x,y)-{\xi}_0(x,y)=y\tilde{\beta}_T(x,y)-y{\beta}_0(x,y)$. Let us stress that the ratio $\frac{\tilde{\pi}_T(X_0,Y_0)}{\pi_0(X_0,Y_0)}$ in the expression of $Z^{(T)}$ comes from the fact that the two diffusions $(X^{(0)},X^{(0)})_t$ and
$(X^{(T)},X^{(T)})_t$ have different initial laws, since they are both stationary, with the different stationary laws.

2) Let us control by below $\mathbb{P}^{(T)}_{V_0,\beta_0}(Z^{(T)} \ge \frac{1}{\lambda})$ for $\lambda >0$. Recalling the definition of $\pi_0$ and $\tilde{\pi}_T$ (see \eqref{E:def_pi_0_beta_0}, \eqref{E:def_tilde_pi_prior})), we see that $\tilde{\pi}_T/\pi_0$ is equal to $1$ outside some compact set (that can be chosen independent of $T$), and converges uniformly to $1$ on this compact set. 
Hence it is bounded away from zero if $T$ is large, and $\frac{\tilde{\pi}_T(X_0,Y_0)}{\pi_0(X_0,Y_0)} \ge C >0$ almost surely.  

Hence, we will focus on the exponential term in \eqref{E:Girsanov_explicit}, that we note $\mathcal{E}^{(T)}=Z^{(T)} 
\frac{\pi_0(X_0,Y_0)}{\tilde{\pi}_T(X_0,Y_0)}$. We know that under $\mathbb{P}^{(T)}_{V_0,\beta_0}$ the canonical process $(X,Y)_t$ has the same law as $(X^{(0)},Y^{(0)})_t$ defined in \eqref{E:def_X_0}--\eqref{E:def_Y_0}. Hence, the law of $\log (\mathcal{E}^{(T)})$ is the law of the random variable
\begin{multline*}
-\frac{1}{4}
\int_0^T \Delta_T(X^{(0)}_s,Y^{(0)}_s) \dd Y^{(0)}_s
\\
-\frac{1}{8 \sigma^2} \int_0^T [\sigma^2 {\xi}_0(X^{(0)}_s,Y^{(0)}_s) + \sigma^2 \Delta_T(X^{(0)}_s,Y^{(0)}_s) + V_0'(X^{(0)}_s,Y^{(0)}_s) ]^2 
\\
-  [\sigma^2 {\xi}_0(X^{(0)}_s,Y^{(0)}_s) + V_0'(X^{(0)}_s,Y^{(0)}_s) ]^2 \dd s.
\end{multline*}
This random variable is equal, using \eqref{E:def_Y_0} and after some computations, to
\begin{align*}
&-\frac{\sigma}{2} \int_0^T \Delta_T(X^{(0)}_s,Y^{(0)}_s) \dd B_s
-\frac{\sigma^2}{8 } \int_0^T \Delta_T(X^{(0)}_s,Y^{(0)}_s)^2\dd s\\
:=
&-M_T-I_T.
\end{align*}
Using the previous considerations we can write that, for $T$ large enough,
\begin{align*}
\mathbb{P}^{(T)}_{V_0,\beta_0}(Z^{(T)} \ge \frac{1}{\lambda})
&\ge \mathbb{P}^{(T)}_{V_0,\beta_0} (\mathcal{E}^{(T)} \ge \frac{1}{C \lambda})
\\
& = \mathbb{P}^{(T)}_{V_0,\beta_0} \left( - \log  \mathcal{E}^{(T)}  \le \log(C \lambda) \right)
\\
&=1-\mathbb{P}^{(T)}_{V_0,\beta_0} \left( - \log  \mathcal{E}^{(T)}   >  \log(C \lambda) \right)
\\
&\ge 1-\mathbb{P}^{(T)}_{V_0,\beta_0} \left( \abs{ \log  \mathcal{E}^{(T)}}   >  \log(C \lambda) \right)
\\
&=1-\mathbb{P}\left( \abs{M_T+I_T} > \log(C \lambda) \right)
\end{align*}
where in the last line we have used 
 that the law of $\log  \mathcal{E}^{(T)} $ under $\mathbb{P}^{(T)}_{V_0,\beta_0} $ is the law of $-M_T-I_T$.
 Assume now that $\lambda>1/C$, then using Markov inequality, we can write 	
\begin{align*}
\mathbb{P} \left( \abs{M_T+I_T} > \log(C \lambda) \right) 
&\le \mathbb{P} \left( \abs{M_T} > \frac{1}{2}\log(C \lambda)
\right)  + \mathbb{P} \left(  \abs{I_T} > \frac{1}{2}\log(C \lambda)  \right) 
\\
&\le \frac{4}{\log(C \lambda)^2} \mathbb{E}(M_T^2)+ \frac{2}{\log(C \lambda)} \mathbb{E}(\abs{I_T})
\end{align*}
Since 	$\mathbb{E}(M_T^2)=2 \mathbb{E}(I_T)$ by Ito's isometry, we see that the condition
\begin{equation}\label{E:bound_on_I_T}
\sup_{T \ge 0} ~\mathbb{E}(I_T) < \infty
\end{equation}
is sufficient to get that there exists  $\lambda_0$ such that for any $T$ large enough we have,
\begin{equation*}
\mathbb{P}^{(T)}_{V_0,\beta_0}(Z^{(T)} \ge \frac{1}{\lambda_0}) \ge 1/2.
\end{equation*}
It remains to check that \eqref{E:bound_on_I_T} holds true.
Recalling that $I_T=\frac{\sigma^2}{8} \int_0^T \Delta_T(X_s^{(0)}, Y_s^{(0)})^2 \dd s$ and using that the process $(X^{(0)}_t,Y^{(0)}_t)_{t \ge 0}$ is stationary, with invariant law $\pi_0$ we have
\begin{equation*}
\mathbb{E}(I_T)=T \frac{\sigma^2}{8} \mathbb{E} \big[ \Delta_T(X_0^{(0)}, Y_0^{(0)})^2 \big]= 
T \frac{\sigma^2}{8} \int_{\mathbb{R}^2} \Delta_T(x,y)^2 {\pi}_0(x,y) \dd x \dd y.
\end{equation*}
Since $\pi_0$ is a bounded function by \eqref{E:def_pi_0_beta_0}, we deduce
\begin{equation*}
\mathbb{E}(I_T) \le C T \int_{\mathbb{R}^2} \Delta_T(x,y)^2  \dd x \dd y.
\end{equation*}
Recalling that by definition $\Delta_T=\tilde{\xi}_T-{\xi}_0$ and using the assumption \eqref{E:cond_Girsanov_lemma_statement} in the statement of the lemma, we deduce that \eqref{E:bound_on_I_T} holds true and the lemma follows.
%
%
%
\end{pf}	
	
\subsection{Proof of Theorem \ref{T:lower_bound_y0_zero}}
The scheme of the proof is similar to the proof of Theorem \ref{T:lower_bound_y0_non_zero}. However, one needs some modifications taking into account that $y_0=0$.
\subsubsection{Constuction of the prior}
The prior is the same as in the proof of Theorem \ref{T:lower_bound_y0_non_zero} except that we need to modify slightly the functions $V_0$ and $h$. Let us give more details. Let $k_1$, $k_2$ and $R >0$. We choose $V_0 : \mathbb{R} \to \mathbb{R} $ a $\mathcal{C}^\infty$ function such that $V_0(x)=x^2$ for $\abs{x}$ large and $V'_0(x)=0$ on a neighbourhood of $x_0$, and we define 
 \begin{equation*}
 \pi_0(x,y)=c_\eta \exp( -\frac{\eta }{2} [ \frac{y^2}{2}+ V_0(x)  ]), \quad  \beta_0(x,y)=\eta, \quad  \xi_0(x,y)=\eta y,
 \end{equation*}
 where {\modar $0<\eta<1/2$} and
 where $c_\eta$ is the constant that make $\pi_0$ a probability measure. 
 The function $\pi_0$ is $\mathcal{C}^\infty$ and it is possible to choose $\eta$ small enough such that
 \begin{equation*}
 \pi_0 \in \mathcal{H}^{k_1,k_2}(R/2).
 \end{equation*}
 We know from Section \ref{Ss:link_drift_statio} that $\pi_0$ is the unique stationary measure of the process $(X^{(0)},Y^{(0)})$ solution to the stochastic differential equation \eqref{E:def_X_0}--\eqref{E:def_Y_0}. If we set $R'=2/\eta$, then recalling Definition \ref{D:class_drift} we have $\beta_0 \in \Sigma^{k_1,k_2}(V,R/2,R'/2)$.
 
Let $h: \mathbb{R}\to \mathbb{R}$ be a $\mathcal{C}^\infty$ function with support on $[-1,1]$ such that,
 \begin{equation} \label{E:pty_h_y0_0}
 h(0)=1,  ~ h'(0)=0, ~ \int_{-1}^1 h(z) \dd z=0, ~ \int_{0}^1 zh(z) \dd z=\int_{-1}^0 zh(z) \dd z=0.
 \end{equation} 
For $T >0$ we define the perturbation of $\pi_0$, as in Section \ref{Ss:Constru_prior_non_zero} by
 \begin{equation*}
 \tilde{\pi}_T(x,y)= \pi_0(x,y) + \frac{1}{M_T} h(\frac{x-x_0}{h_1(T)}) h(\frac{y}{h_2(T)}),
 \end{equation*}
 where  $M_T\to \infty$, $h_1(T)\to 0$, $h_2(T) \to 0$ will be calibrated latter. Again $\tilde{\pi}_T$ is a smooth probability measure for $T$ large enough and we define 
 \begin{equation*}
 \tilde{\beta}_T(x,y)= \beta_{\tilde{\pi}_T}(x,y), \quad  \tilde{\xi}_T(x,y) = y \tilde{\beta}_T(x,y)=\xi_{\tilde{\pi}_T}(x,y),
 \end{equation*}
 where we used the definitions \eqref{E:def_xi} and \eqref{E:def_beta_g}.

The following lemma gives an assessment of the difference between $\beta_0$ and $\tilde{\beta}_T$. 
 \begin{lem}\label{L:beta_beta_T_y_0_zero}
 	1) Recall the definition of the following compact set of $\mathbb{R}^2$
 	\begin{equation*}
 	K_T=[x_0-h_1(T),x_0+h_1(T)]\times [-h_2(T),h_2(T)].
 	\end{equation*}
 	Then, for $T$ large enough, we have for all $(x,y) \notin K_T$ :
 	\begin{equation*}
 	\beta_0(x,y)=\tilde{\beta}_T(x,y),\quad \xi_0(x,y)=\tilde{\xi}_T(x,y).
 	\end{equation*}
 	2) For $(x,y) \in K_T$, we have the control,
 	\begin{align}
 	\label{E:maj_diff_beta_y0_zero}
 	&\abs{ \beta_0(x,y)-\tilde{\beta}_T(x,y)}  \le \frac{C}{M_T}\left[\frac{h_2(T)}{h_1(T)}+ \frac{1}{h_2(T)^2}\right],
 	\\ \label{E:maj_diff_xi_y0_zero}
 	&\abs{\xi_0(x,y)-\tilde{\xi}_T(x,y)}  \le \frac{C}{M_T}\left[\frac{h_2(T)^2}{h_1(T)}+ \frac{1}{h_2(T)}\right],
 	\end{align}
 	where $C$ is some constant independent of $T$, $h_1(T)$, $h_2(T)$, $M_T$.
 	
 	3) We have
 	\begin{equation*}
 	\int_{\mathbb{R}^2} \abs{ \tilde{\xi}_T(x,y) - \xi_0(x,y)     }^2 \dd x \dd y
 	\le 
 	\frac{C}{M_T^2} \left[\frac{h_2(T)^5}{h_1(T)}+ \frac{h_1(T)}{h_2(T)}\right].
 	\end{equation*}
 \end{lem}
 \begin{pf}
 	1) We first prove that $\tilde{\xi}_T$ and $\xi_0$ coincides on $K_T^c$. Using the notations and arguments of Lemma \ref{L:beta_beta_T}, we know that $\tilde{\xi}_T(x,y)=\xi_0(x,y),$ for all $(x,y) \notin K_T$ is a consequence of $\mathcal{I}[d_T](x,y)=0$ for $(x,y) \notin K_T$. We recall that $\mathcal{I}=\sum_{i=1}^3 \mathcal{I}^i$ is given by \eqref{E:def_I_1_g}--\eqref{E:def_I_3_g} and $d_T(x,y)=\tilde{\pi}_T(x,y) - \pi_0(x,y)$ is given by :
 	\begin{equation} \label{E:express_d_T}
 	d_T(x,y)=\frac{1}{M_T}h(\frac{x-x_0}{h_1(T)})h(\frac{y}{h_2(T)}).
 	\end{equation}
 	If $(x,y) \notin K_T$, the first situation is $\abs{x-x_0}>h_1(T)$, then $d_T(x,z)=\frac{\partial d_T}{\partial x}(x,z)=\frac{\partial d_T}{\partial y}(x,z)=0$ for all $z \in \mathbb{R}$ and we deduce that $\mathcal{I}^i[d_T](x,y)=0$ for $i=1,2,3$. The second situation is $\abs{y}>h_2(T)$ and $\abs{x-x_0}\le h_1(T)$. In that case $\mathcal{I}^2[d_T](x,y)=0$ for $T$ large enough, by using that from assumption on $V_0$, $V'_0(x)=0$ for $x$ in some neighbourhood of $x_0$. From \eqref{E:def_I_3_g}, we have $\mathcal{I}^3[d_T](x,y)=\frac{2}{M_T h_2(T)} h(\frac{x-x_0}{h_1(T)})[h'(\frac{y}{h_2(T)})-h'(0)]$ which is equal to $0$ since $\abs{y}>h_2(T)$ and $h'(0)=0$ by \eqref{E:pty_h_y0_0}. In order to check that $\mathcal{I}^1[d_T](x,y)=0$, let us assume for simplicity that $y>h_2(T)$, as the case $y<-h_2(T)$ is similar. Then,
 	\begin{multline*}
 	\mathcal{I}^1[d_T](x,y)=\frac{1}{\sigma^2 h_1(T)} \int_0^{h_2} z h(\frac{z}{h_2(T)}) \dd z
 	 h'(\frac{x-x_0}{h_1(T)}) \\
 	 =\frac{h_2(T)^2}{\sigma^2 h_1(T)} \int_0^1 z h(z) \dd z h'(\frac{x-x_0}{h_1(T)})=0
 	\end{multline*}
 	by \eqref{E:pty_h_y0_0}. 
Eventually, this gives that	$\mathcal{I}[d_T](x,y)=0$ for $(x,y) \notin K_T$, and thus $\tilde{\xi}_T(x,y)=\xi_0(x,y)$. 

The equality between the functions $\tilde{\beta}_T$ and $\beta_0$ on $K_T^c$ is a consequence of $\tilde{\beta}_T(x,y)=\tilde{\xi}(x,y)/y$, $\beta_0(x,y)={\xi_0}(x,y)/y$ for $y \neq 0$.

2) We first prove \eqref{E:maj_diff_xi_y0_zero}. Recalling \eqref{E:diff_xi_tilde_xi_0}, the fact that $\norme{\pi_0-\tilde{\pi}}_\infty=\norme{d_T}_\infty \le C/M_T$, and that $\tilde{\pi}_T$ is lower bounded on $K_T$ as soon as $T$ is large enough, we deduce
 \begin{align}\nonumber
 \forall (x,y) \in K_T,~ \abs{\tilde{\xi}_T(x,y)-\xi_0(x,y)}  & \le C \big[\frac{\abs{\xi_0(x,y)}}{M_T} + \mathcal{I}[d_T](x,y) \big]
 \\ \label{E:maj_diff_xi_intermediate_y0_zero}
 &\le C \big[\frac{\abs{y}}{M_T} + \mathcal{I}[d_T](x,y) \big].
 \end{align}
Now, we use $\abs{\mathcal{I}[d_T](x,y) } \le \abs{ \mathcal{I}^1[d_T](x,y)  } +  \abs{ \mathcal{I}^2[d_T](x,y)  } +  \abs{ \mathcal{I}^3[d_T](x,y)  }$. From \eqref{E:def_I_1_g}, we have $\abs{\mathcal{I}^1[d_T](x,y) }  \le C \abs{\int_0^y z \dd z} 
\norme{ \frac{\partial d_T}{ \partial x} }_\infty \le C \abs{y} h_2(T) \frac{C}{M_T h_1(T)}$, for all $(x,y) \in K_T = [x_0-h_1(T),x_0+h_1(T)]\times [-h_2(T),h_2(T)]$, and where we have used the expression \eqref{E:express_d_T} for $d_T$.
As $x \mapsto V'(x)$ vanishes on a neighbourhood of $x_0$, we get that for $T$ large enough $\mathcal{I}^2[d_T](x,y)=0$ for $(x,y) \in K_T$. From \eqref{E:def_I_3_g}, we deduce $\abs{ \mathcal{I}^3[d_T](x,y)} \le C \abs{y} \norme{\frac{\partial^2 d_T}{\partial y^2} }_\infty \le 
 C  \frac{\abs{y}}{M_T h_2(T)^2} $. Collecting the controls on $\mathcal{I}^i[d_T](x,y)$ for $i=1,2,3$, with \eqref{E:maj_diff_xi_intermediate_y0_zero} we get
 \begin{equation*}
  \forall (x,y) \in K_T,~ \abs{\tilde{\xi}_T(x,y)-\xi_0(x,y)} \le C \frac{\abs{y}}{M_T} \big[ 1 + \frac{h_2(T)}{h_1(T)} + \frac{1}{h_2(T)^2}  \big].  
 \end{equation*}
Using that for $(x,y) \in K_T$, we have $\abs{y} \le h_2(T)$ and the last equation implies \eqref{E:maj_diff_xi_y0_zero}. Moreover, from the fact that $\tilde{\xi}_T(x,y)-\xi_0(x,y)=y [ \tilde{\beta}_T(x,y)-\beta_0(x,y)]$, it implies \eqref{E:maj_diff_beta_y0_zero}.

3) The  third point of the lemma in a consequence of the first two points and the fact that the Lebesgue measure of $K_T$ is proportional to  $h_1(T)h_2(T)$. 
\end{pf}

We now state a result analogous to Lemma \ref{L:beta_in_prior}, but in the situation $y_0=0$.  
\begin{lem} \label{L:beta_in_prior_y_0_zero}
	Let $\varepsilon >0$ and assume that for all $T$ large,
	\begin{equation*} 
	M_T^{-1} \le \varepsilon h_1(T)^{k_1}, \quad M_T^{-1} \le \varepsilon h_2(T)^{k_2},
	\end{equation*}
	and
	\begin{equation*}
	\frac{h_2(T)}{h_1(T)} + \frac{1}{h_2(T)^2} = o(M_T) \text{ as $T \to \infty$.}
	\end{equation*}
	Then, if $\varepsilon>0$ is small enough, 
	we have 
	\begin{equation*}
	\tilde{\beta}_T \in \Sigma^{k_1,k_2}(V_0,R,R'), 
	\end{equation*} 
	for all $T$ sufficiently large.	
\end{lem}
We omit the proof of Lemma \ref{L:beta_in_prior_y_0_zero} as it is similar to the proof of Lemma \ref{L:beta_in_prior} (except that we use \eqref{E:maj_diff_beta_y0_zero} instead of \eqref{E:maj_diff_beta}).

\subsubsection{Proof of the lower bound \eqref{E:lower_bound_y0_zero} on the minimax risk}
We omit most of the details of the proof as it is very similar to the proof given in Section \ref{Ss:proof_lower_bound_y_0_non_zero}.
Indeed, by  repeating the arguments of the proof given in Section \ref{Ss:proof_lower_bound_y_0_non_zero}, relying on Lemmas \ref{L:beta_beta_T_y_0_zero}--\ref{L:beta_in_prior_y_0_zero} instead of Lemmas \ref{L:beta_beta_T}--\ref{L:beta_in_prior}, we deduce that,
\begin{equation} \label{E:generic_lower_bound_y_0_non_zero}
R(\widetilde{\pi}_T(x_0,0)) \ge \frac{C}{M_T^2},
\end{equation}
as soon as we can find $\varepsilon>0$, $h_1(T)\to 0$, $h_2(T)\to 0$ and $M_T\to \infty$ satisfying the conditions
\begin{align} \label{E:cond1_y_0_zero}
&	M_T^{-1} \le \varepsilon h_1(T)^{k_1}, \quad M_T^{-1} \le \varepsilon h_2(T)^{k_2},
	\\
	\label{E:cond2_y_0_zero}
&		\frac{h_2(T)}{h_1(T)} + \frac{1}{h_2(T)^2} = o(M_T) \text{ as $T \to \infty$},
	\\
	\label{E:cond3_y_0_zero}
&		\sup_T \frac{T}{M_T^2}\left[\frac{h_2(T)^5}{h_1(T)}+ \frac{h_1(T)}{h_2(T)}\right] < \infty.
\end{align}
Let us maximise $1/M_T^2$ under these three constraints.

\underline{Case 1, $k_1 < k_2/3$ :} we set $h_1(T)=h_2(T)^3$ and $h_2(T)=\left( \frac{1}{\varepsilon M_T} \right)^{1/k_2}$ and the conditions \eqref{E:cond1_y_0_zero} hold true, using $k_1<k_2/3$. With these choices, the condition \eqref{E:cond3_y_0_zero} reduces to the boundedness of $\frac{T}{M_T^2} M_T^{-2/k_2}$, which is implied if we set $M_T=T^{-\frac{k_2}{2+2k_2}}$. 
Next \eqref{E:cond2_y_0_zero} holds true as $ \frac{h_2(T)}{h_1(T)} + \frac{1}{h_2(T)^2}  = O( \frac{1}{h_2(T)^2}) = O( M_T^{2/k_2})=o(M_T) $ as $k_2 >2$. Hence, we can
use \eqref{E:generic_lower_bound_y_0_non_zero} with $M_T=T^{-\frac{1}{2+2k_2}}$, yielding to
\eqref{E:lower_bound_y0_zero} in the case $k_1 < k_2/3$.

\underline{Case 2, $k_1 \ge  k_2/3$ :}  we set $h_1(T)=h_2(T)^3$ and $h_1(T)=\left( \frac{1}{\varepsilon M_T} \right)^{1/k_1}$ and the conditions \eqref{E:cond1_y_0_zero} follows. Now, the condition \eqref{E:cond3_y_0_zero}  reduces to the boundedness 
of $\frac{T}{M_T^2} M_T^{-2/(3k_1)}$, which yields to the choice $M_T=T^{-\frac{k_1}{2k_1+2/3}}$.
Next \eqref{E:cond2_y_0_zero} holds true as $ \frac{h_2(T)}{h_1(T)} + \frac{1}{h_2(T)^2}  = O( \frac{1}{h_2(T)^2}) = O( M_T^{2/(3k_1)})=o(M_T) $ as $k_1 >2/3$. Eventually, we deduce \eqref{E:lower_bound_y0_zero} from \eqref{E:generic_lower_bound_y_0_non_zero}.

%
 
%
%


\section{Appendix} \label{Appendix}

{\modar In this section we prove the technical Lemmas \ref{L:HShort}--\ref{L:HGlobal} on the semi group of the process.
}
	
\subsection{Proof of Lemma \ref{L:HShort}}
{\modar 
This proof is exactly the same as the one of Corollary 2.12 in \cite{Cattiaux_et_al_14}, after remarking 
that the results of Theorem 2.1. in \cite{Konakov_et_al_10} can be applied
to S.D.E. with $\mathcal{C}^1$ coefficients.
\qed
}

{\modar 
\subsection{Proof of Lemma \ref{L:HGlobal}}
	We first prove that \eqref{E:equa_type_HGlob} holds for $t=D$. Let us denote $\widetilde{K}=\{z \in \mathbb{R}^2 \mid d(z,K) \le 1 \}$ the compact set of points at distance less than $1$ of $K$. Since $D<1$, we can apply Lemma \ref{L:HShort} with the choice of compact set $\widetilde{K}$, to get that if $f$ has support on $K \subset \widetilde{K}$, and $z \in \widetilde{K}$
	\begin{align} \nonumber
	\abs{P_D(f)(z)} & \le \int_{\mathbb{R}^2} \abs{f(z')} p_D(z;z') dz'
	\\ \nonumber
	 &\le  \int_{\mathbb{R}^2} \abs{f(z')} p_D^G(z;z') dz' +   \int_{\mathbb{R}^2} 
	\abs{f(z')} p_D^U(z;z') dz'
	\\ \label{E:proofHGlobal_controlP_z_small}
	& \le \frac{C_G}{D^2} \int_{\mathbb{R}^2} \abs{f(z')} dz' + C_U \norme{f}_\infty e^{-1/(C_U D)}.
	\end{align}
	Hence, it proves \eqref{E:equa_type_HGlob} for $t=D$ and $z \in \widetilde{K}$. 
	
	If $z \notin \widetilde{K}$, we let $T_{\widetilde{K}}=\inf \{t \ge 0 \mid Z_t \in \widetilde{K}\}$ the entrance time in the compact set $\widetilde{K}$, which is a stopping time.
	As the support of $f$ is included in $K$, we have by continuity of the process, $P_D(f)(z)=E_z[f(Z_D)]=E_z[f(Z_D)1_{\{T_{\widetilde{K}} \le t\}}]$. Using the strong Markov property at the time $T_{\widetilde{K}}$ we deduce, 
	\begin{equation} \label{E:proofHGlobalMarkov}
	P_D(f)(z)=E_z[P_{D-T_{\widetilde{K}}}(f)(Z_{T_{\widetilde{K}}}) 1_{\{T_{\widetilde{K}} \le D\}}].
	\end{equation}
	By the continuity of the process, we remark that $d(Z_{T_{\widetilde{K}}},K)=1$ on the set $T_{\widetilde{K}}\le D$, and $D-T_{\widetilde{K}}$ in $(0,D)$ as $z \notin \widetilde{K}$. This lead us to consider for $z'\in \widetilde{K}$ with $d(z',K)=1$ and $s \in (0,D) \subset (0,1)$, an upper bound for
	\begin{align} \nonumber
	\abs{P_s(f)(z')} &\le \int_{\mathbb{R}^2} \abs{f(w)} p_s^G(z';w) dw +   \int_{\mathbb{R}^2} 
	\abs{f(w)} p_s^U(z';w) dw
	\\ \label{E:proofHGlobalPzprime}
	& \le \int_{\mathbb{R}^2}\abs{f(w)} dw \tilde{C}_G + C_U \norme{f}_\infty e^{-1/(C_U s)},
	\end{align}
	where $\tilde{C}_G = \sup \{ p_s^G(z',w) \mid s \in (0,1), w \in K, z' \in \widetilde{K} \text{ with } d(z',K)=1 \}$, and where we used again Lemma \ref{L:HShort}.	We can see that $\tilde{C}_G$ is finite. Indeed, if $z'=(x',y')$ is such that  $d(z',K)=1$ and $w=(w_1,w_2) \in K$, we have
	$p_s^G(z',w) \le \frac{C_G}{s^2} \exp (-\frac{1}{C_G} [\frac{(w_2-y')^2}{s}+\frac{(w_1-x'-\frac{w_2+y'}{2}s)^2}{s^3} ])$. Using the inequality 
	$A^2 \le (A-B)^2(1+1/s)+B^2(1+s)$ for any $A$, $B$, that entails
	$(A-B)^2 \ge A^2 \frac{s}{s+1} - B^2 s$,	 we deduce
	$\frac{(w_2-y')^2}{s}+\frac{(w_1-x'-\frac{w_2+y'}{2}s)^2}{s^3}  \ge 
	\frac{(w_2-y')^2}{s} + \frac{(w_1-x')^2}{s^3} \frac{s}{s+1} - \frac{(w_2+y')^2}{2s} s$. Using that $s<1$ and that $\abs{w_2}$ and $\abs{y'}$ are bounded by some constant depending on the compact $K$, we deduce $\frac{(w_2-y')^2}{s}+\frac{(w_1-x'-\frac{w_2+y'}{2}s)^2}{s^3}  \ge d(z',w)^2/(2s) - \tilde{C}_K$ for some constant $\tilde{C}_K$ depending on the compact $K$ only. It gives $p_s^G(z',w) \le \frac{C_G}{s^2} \exp (-\frac{d(z',w)^2}{C_G 2 s}  ) \exp( \frac{\widetilde{C}_K}{C_G} )$. As $w\in K$ and $d(w',K)=1$, we deduce that $p_s^G(z',w) \le \frac{C_G}{s^2} \exp (-\frac{1}{C_G 2 s}  ) \exp( \frac{\widetilde{C}_K}{C_G} )$ and thus $\widetilde{C}_G$ is finite.
	Joining \eqref{E:proofHGlobalPzprime} and \eqref{E:proofHGlobalMarkov}, we deduce that, for $z \notin \widetilde{K}$
	\begin{align} \nonumber
	\abs{P_D(f)(z)} &\le \widetilde{C}_G \norme{f}_{L^1(\mathbb{R}^2)} +  C_U \norme{f}_\infty E_z[ e^{-1/(C_U(D-T_{\widetilde{K}}))} 1_{\{T_{\widetilde{K}} \le D\}}]
	\\ \label{E:proofHGlobal_controlP_z_large}
	&\le \widetilde{C}_G \norme{f}_{L^1(\mathbb{R}^2)} +  C_U \norme{f}_\infty e^{-1/(C_U D)}.
	\end{align}
	The control  \eqref{E:equa_type_HGlob} for $t=D$ is now a consequence of \eqref{E:proofHGlobal_controlP_z_small} and \eqref{E:proofHGlobal_controlP_z_large}.
	Eventually, we prove that \eqref{E:equa_type_HGlob} for $t=D$ is sufficient to deduce the lemma. Let $0<D<1$ and $t>D$, then for $z \in \mathbb{R}^2$, we write
	\begin{equation*}
	\abs{P_t(f)(z)} \le \int_{\mathbb{R}^2} p_{t-D}(z,z') \abs{P_D(f)(z')} dz'
	\end{equation*}
	and using the estimate \eqref{E:equa_type_HGlob} for $\abs{P_D(f)(z')}$ gives the result for $\abs{P_t(f)(z)}$. 
	\qed
}

\end{document}